\documentclass[11pt]{amsart}
\usepackage[top=1in, bottom=0.9in, left=1in, right=1in]{geometry}

\usepackage{mathtools,mathrsfs,amsthm,amsfonts,amssymb,graphicx,appendix,cancel}
\usepackage{xspace,enumerate,comment,bbm,nth}

\usepackage{tikz,pgfplots}

\usepackage{hyperref}
\usepackage[normalem]{ulem}

\usepackage[active]{srcltx}

\numberwithin{equation}{section}
\newtheorem{theorem}{\sc Theorem}[section]

\newtheorem{corollary}[theorem]{\sc Corollary}
\newtheorem{definition}[theorem]{\sc Definition}
\newtheorem{lemma}[theorem]{\sc Lemma}

\theoremstyle{plain}

\newtheorem{prop}[theorem]{\sc Proposition}
\newtheorem{cor}[theorem]{\sc Corollary}

\theoremstyle{remark}
\newtheorem{remark}[theorem]{\sc Remark}

 \newcommand{\ol}[1]{\overline{#1}}
 \newcommand{\ul}[1]{\underline{#1}}
 \newcommand{\wh}[1]{\widehat{#1}}

\newcommand{\F}{\mathcal F}
\newcommand{\R}{\mathbb{R}}
\newcommand{\Sol}{\mathcal S}
\newcommand{\Z}{\mathbb{Z}}
\newcommand{\N}{\mathbb{N}}

\renewcommand{\P}{\mathbb{P}}
\newcommand{\E}{\mathbb{E}}

\newcommand{\cal}[1]{\mathcal #1}

\newcommand{\HV}{{\mathcal H}}
\newcommand{\HF}{{\mathcal H}}

\newcommand{\eps}{\varepsilon}
\renewcommand{\epsilon}{\varepsilon}

\renewcommand{\emptyset}{\varnothing}

\renewcommand{\le}{\leqslant}
\renewcommand{\ge}{\geqslant}

\newcommand{\CC}{\D C}

\newcommand{\cyl}{(0,+\infty)\times\R}
\newcommand{\ccyl}{[0,+\infty)\times\R}

\newcommand{\Ham}{\mathscr{G}}
\newcommand{\Hamone}{\mathscr{G}_1}
\newcommand{\Hamzero}{\mathscr{G}_0}
\newcommand{\HamzeroN}{\mathscr{G}_{0,N}}
\newcommand{\Fam}{\mathscr{F}}

\newcommand{\D}[1]{\mbox{\rm #1}}

\newcommand{\essinf}{\mathrm{ess\,inf\,}}
\newcommand{\esssup}{\mathrm{ess\,sup\,}}

\newcommand{\ucv}{\rightrightarrows}

\begin{document}
	
\title[Homogenization of 1-D viscous HJ
equations]{Stochastic homogenization of nonconvex viscous
  Hamilton-Jacobi equations in one space dimension}

\author[A.\ Davini]{Andrea Davini}
\address{Andrea Davini\\ Dipartimento di Matematica\\ {Sapienza} Universit\`a di
  Roma\\ P.le Aldo Moro 2, 00185 Roma\\ Italy}
\email{davini@mat.uniroma1.it}
\urladdr{https://www1.mat.uniroma1.it/people/davini/home.html}

\author[E.\ Kosygina]{Elena Kosygina} \address{Elena Kosygina\\
  Department of Mathematics\\ Baruch College\\ One Bernard Baruch
  Way\\ Box B6-230, New York, NY 10010\\ USA}
\email{elena.kosygina@baruch.cuny.edu}
\urladdr{http://www.baruch.cuny.edu/math/elenak/} \thanks{E.\ Kosygina
  was partially supported by the Simons Foundation (Award \#523625)
  and by the Visiting Professor Programme 2020 at Sapienza
  Universit\`a di Roma.}

\author[A.\ Yilmaz]{Atilla Yilmaz}
\address{Atilla Yilmaz\\ Department of Mathematics\\ Temple University\\ 1805 North Broad Street, Philadelphia, PA 19122, USA}
\email{atilla.yilmaz@temple.edu}
\urladdr{http://math.temple.edu/$\sim$atilla/}
\thanks{A.\ Yilmaz was partially supported by the Simons Foundation (Award \#949877).}

\date{\today}


\begin{abstract}
  We prove homogenization for viscous Hamilton-Jacobi equations with a
  Hamiltonian of the form $G(p)+V(x,\omega)$ for a wide class of
  stationary ergodic random media in one space dimension.  The
  momentum part $G(p)$ of the Hamiltonian is a general (nonconvex)
  continuous function with superlinear growth at infinity, and the
  potential $V(x,\omega)$ is bounded and Lipschitz continuous. The
  class of random media we consider is defined by an explicit hill
  and valley condition on the diffusivity-potential pair which is
  fulfilled as long as the environment is not ``rigid''.
\end{abstract}

\maketitle

\section{Introduction}\label{sec:intro}

We are interested in the homogenization as $\epsilon\to 0+$ of viscous
Hamilton-Jacobi (HJ) equations
\begin{equation}\label{eq:introHJ}
  \partial_t u^\epsilon=\epsilon a\left(\frac{x}{\epsilon},\omega\right) \partial^2_{xx} u ^\epsilon+G(\partial_x u^\epsilon)+  \beta V\left(\frac{x}{\epsilon},\omega\right),  \quad (t,x)\in (0,+\infty)\times\R,
\end{equation}
where $(a(x,\omega))_{x\in\R}$ and $(V(x,\omega))_{x\in\R}$ are
realizations of continuous stationary ergodic processes defined on
some probability space $(\Omega,{\cal F},\P)$, the function $G:\R\to\R$ is
continuous and has superlinear growth at $\pm\infty$, and $\beta$ is a
positive constant adjusting the magnitude of $V$. The full set of
assumptions is deferred to Section \ref{sec:result}, but the emphasis
of this work is on nonconvex $G$.

Equations of the form \eqref{eq:introHJ} are a subclass of general
viscous stochastic HJ equations
\begin{equation}\label{vHJ}
  \partial_tu^\epsilon= \epsilon\text{tr}\left(a\left(\frac{x}{\epsilon},\omega\right) D^2 u^\epsilon\right)+H\left(D u^\epsilon,\frac{x}{\epsilon},\omega\right),  \quad (t,x)\in (0,+\infty)\times\R^d,
\end{equation}
where $a(\cdot,\omega)$ is a bounded, symmetric, and nonnegative
definite $d\times d$ matrix with a Lipschitz square root (uniformly in
$\omega\in\Omega$).  Solutions to all PDEs throughout the paper will
be understood in the viscosity sense (see,
\cite{users,barles_book,bardi}).

We shall say that \eqref{vHJ} {\em homogenizes} if there exists a continuous
nonrandom function $\ol{H}:\R^d\to\R$, the effective
  Hamiltonian, such that, with probability 1, for every uniformly
continuous function $g$ on $\R^d$, the solutions $u^\epsilon$ of
\eqref{vHJ} satisfying $u^\epsilon(0,\,\cdot\,,\omega) = g$ converge
locally uniformly on $[0, +\infty)\times \R^d$ as $\epsilon\to 0+$ to
the unique solution $\ol{u}$ of the equation
\begin{equation}\label{effeq}
  \partial_t \ol{u}=\ol{H}(D\ol{u}), \quad (t,x)\in (0,+\infty)\times\R^d,
\end{equation}
satisfying $\ol{u}(0,\,\cdot\,) = g$.

The main result of our paper, Theorem~\ref{thm:genhom}, asserts that, under a mild additional assumption on the pair of processes $(a(x,\omega),V(x,\omega))_{x\in\R}$, equation \eqref{eq:introHJ} homogenizes for a general continuous $G(p)$ with superlinear growth as $|p|\to\infty$. The precise statement is given in Section~\ref{sec:result}.

The literature on the homogenization of viscous and inviscid
($a\equiv 0$) HJ equations in stationary ergodic media is vast, and we
refer, for instance, to \cite{AT15,KYZ20,DK22} for a broader review of
the existing body of work on what is often called qualitative
homogenization, i.e., results which show the existence of the limit
as $\epsilon\to 0+$ but do not quantify the rate of
convergence. Quantitative results, in particular algebraic rates of
convergence for \eqref{vHJ} under various sets of assumptions, were
obtained in \cite{ACS14,AC15,AC18}.

Our interest in the one-dimensional setting is motivated by two
considerations. First, it was shown by counterexamples
(\cite{Zil,FS,FFZ}) that, for dimensions $d\ge 2$, homogenization can
fail if the Hamiltonian has a strict saddle point (otherwise being
``standard'', in particular superlinear as $|p|\to+\infty$ uniformly
in $(x,\omega)$). Therefore, unlike in the periodic case, there can be
no general homogenization result for nonconvex superlinear
Hamiltonians in the stationary ergodic setting for $d\ge 2$, at least,
not without additional assumptions on the mixing properties of the
environment or the ``shape" of the Hamiltonian (e.g., rotational
symmetry or homogeneity in the momentum variables). Second, for $d=1$
and $a\equiv 0$, the general HJ equation with a continuous coercive
Hamiltonian homogenizes (\cite{ATY_1d,Gao16,Y21a}), and a similar
result is also expected for $a\not\equiv 0$. The current paper
establishes such a result for an equation of the form
\eqref{eq:introHJ} in stationary ergodic environments which satisfy
the scaled hill and valley condition on $(a,V)$ (see Section
\ref{sec:result}). We note that this last condition is automatically
satisfied when $a\equiv 0$ (see Remark \ref{rem:vac}). Therefore, our
result can be considered as a viscous version of the one in
\cite{ATY_1d}.

The viscous case, even in the uniformly elliptic setting, is known to
present additional challenges which cannot be overcome by mere
modifications of the methods used for $a\equiv 0$. A representative
example of the above is \eqref{vHJ} with a Hamiltonian that is
quasiconvex in $p$\footnote{i.e., $\forall x,\omega,\lambda$, the
  sublevel set $\{p\in\R^d:\ H(p,x,\omega)\le \lambda\}$ is convex.}:
while homogenization in the inviscid case has been known for quite
some time for all $d\ge 1$ (\cite{DS09}, \cite{AS}), the corresponding
result in the viscous case has only been proven for $d=1$
(\cite{Y21b}, $a>0$) and it remains an open problem for $d\ge 2$.

Returning to stochastic homogenization of viscous HJ equations in one
space dimension, we shall briefly review currently known homogenizable classes
of Hamiltonians and point out remaining challenges. In \cite[Section
4]{DK17}, the first two authors proved the homogenization of
\eqref{vHJ} ($d=1$) with $H(p,x,\omega)$ which are ``pinned'' at
finitely many points $p_1<p_2<\ldots<p_n$ (i.e.,
$H(p_i,x,\omega)\equiv \mathrm{const}_i$, $i\in\{1,2,\ldots,n\}$) and
convex on each interval $(-\infty,p_1)$, $(p_1,p_2)$, \ldots,
$(p_n,+\infty)$. For example, the Hamiltonian
\[H(p,x,\omega)=|p|^\gamma-c(x,\omega)|p| =\min\{|p|^\gamma-c(x,\omega)p,|p|^\gamma+c(x,\omega)p\}\] is pinned at a single point $p=0$.
Addition of a nonconstant potential to the equation removes
pinning. In particular, the homogenization of equation \eqref{vHJ}, where
$d=1$, $a(x,\omega)\equiv \mathrm{const} > 0$, and
\begin{equation}
  \label{open}
  H(p,x,\omega):=\frac12\,|p|^2-c(x,\omega)|p|+\beta V(x,\omega),\quad
  0<c(x,\omega)\le C,\quad \text{and }\ \beta>0,
\end{equation}
remained an open problem even when $c(x,\omega)\equiv c>0$. In \cite{YZ19}, the third
author and O.\;Zeitouni were able to handle the case
$c(x,\omega)\equiv c>0$ in the discrete setting of controlled random
walks in a random potential on $\Z$ under an additional assumption on
$V$ which they called the hill and valley condition. This work paved the way for \cite{KYZ20} which proved
homogenization for \eqref{open} with $c(x,\omega)\equiv c$
(and $a(x,\omega)\equiv 1/2$) under the hill and valley condition. While the last
condition is fulfilled for a wide class of typical random
environments without any restriction on their mixing properties, it is
not satisfied if the potential is ``rigid'', for example, in
the periodic case\footnote{See \cite[Appendix B]{DK22} for a
  detailed discussion and examples.}.

The approach of \cite{YZ19,KYZ20} relies on the Hopf-Cole
transformation, stochastic control representations of solutions, and
the Feynman-Kac formula. It is applicable only to \eqref{eq:introHJ}
with $G(p)=\frac12|p|^2-c|p|=\min\{\frac12p^2-cp,\frac12p^2+cp\}$. The
first two authors (\cite{DK22}) found a different proof which uses
only PDE methods. This new approach allows general $a(x,\omega)\ge 0$ and
any $G(p)$ which is a minimum of a finite number of convex
superlinear functions $G_i$ as long as all of them have the same
minimum. The convexity assumption on $G_i$ was used to ensure that
\eqref{eq:introHJ} with $G_i$ in place of $G$ is homogenizable for
each $i$. As in \cite{KYZ20}, the shape of the effective Hamiltonian associated with 
$G$ is derived from the ones associated to each $G_i$. 
This is where the fact that the functions $G_i$ have the 
same minimum plays a critical role. In \cite{Y21b}, the third author has shown homogenization
when $G(p)$ is a quasiconvex superlinear function and $a(x,\omega)>0$.

The main novelty of the current paper is in allowing a general
superlinear $G$ without any restrictions on its ``shape'' or the
number of its local extrema. Analogously to \cite{KYZ20, DK22}, the
function $G$ can be seen as a minimum of quasiconvex superlinear
functions $G_i$, but when these latter have distinct minima, there is
a nontrivial interaction among them in the homogenization process, and
the shape of the effective Hamiltonian associated with $G$ can no
longer be guessed from those associated with each $G_i$. Due to the
presence of the diffusive term, dealing with this interaction as well
as ``gluing'' together separate pieces of the effective Hamiltonian is
much more delicate than in the inviscid case \cite{ATY_1d,Gao16}.
This is where we need the additional assumption on the
environment. In \cite{Y21b} and \cite{DK22}, the
original hill and valley condition on $V$ has been weakened to the
scaled hill and valley condition on $(a,V)$, which is also retained in
the current paper.   Removing this condition and extending the result to
general non-separated Hamiltonians (i.e., not necessarily of the form
$H(p,x,\omega)=G(p)+\beta V(x,\omega)$) remain open problems for
$d=1$.

\section{Main result and an overview of the proof}

\subsection{Assumptions and the main result}\label{sec:result}

We will denote by $\CC(\R)$ and ${\CC}^1(\R)$ the family of continuous functions and functions of class $\CC^1$ on $\R$, respectively. We will regard $\CC(\R)$ as a Polish space endowed with a metric inducing the topology of uniform convergence on compact subsets of $\R$. We will add the subscript $b$ to indicate that we are considering functions that are also bounded on $\R$.

The triple $(\Omega,\F, \P)$ denotes a probability space, where $\Omega$ is a Polish space, ${\cal F}$ is the $\sigma$-algebra of Borel subsets of $\Omega$, and $\P$ is a complete probability measure on $(\Omega,{\cal F})$.\footnote{The assumptions that $\Omega$ is a Polish space and $\P$ is a complete probability measure are used only in the proof of Lemma \ref{lem:measurable}.} We will denote by ${\cal B}$ the Borel $\sigma$-algebra on $\R$ and equip the product space $\R\times \Omega$ with the product $\sigma$-algebra ${\mathcal B}\otimes {\cal F}$.

We will assume that $\P$ is invariant under the action of a one-parameter group $(\tau_x)_{x\in\R}$ of transformations $\tau_x:\Omega\to\Omega$. More precisely, we assume that the mapping
$(x,\omega)\mapsto \tau_x\omega$ from $\R\times \Omega$ to $\Omega$ is measurable, $\tau_0=id$, $\tau_{x+y}=\tau_x\circ\tau_y$ for every $x,y\in\R$, and $\P\big(\tau_x (E)\big)=\P(E)$ for every $E\in{\cal F}$ and $x\in\R$. We will assume in addition that the action of $(\tau_x)_{x\in\R}$ is {\em ergodic}, i.e., any measurable function $\varphi:\Omega\to\R$ satisfying $\P(\varphi(\tau_x\omega) = \varphi(\omega)) = 1$ for every fixed $x\in{\R}$ is almost surely equal to a constant.

A random process $f:\R\times \Omega\to \R$ is said to be {\em stationary} with respect to $(\tau_x)_{x\in\R}$ if  $f(x,\omega)=f(0,\tau_x\omega)$ for all $(x,\omega)\in\R\times\Omega$. Moreover, whenever the action of $(\tau_x)_{x\in\R}$ is ergodic, we refer to $f$ as a stationary ergodic process.

In this paper, we will consider an equation of the form 
\begin{equation}\label{eq:generalHJ}
\partial_t u=a(x,\omega) \partial^2_{xx} u +G(\partial_x u)+  \beta V(x,\omega),\quad (t,x)\in (0,+\infty)\times\R,
\end{equation}
where $\beta\geqslant 0$, and $a:\R\times\Omega\to (0,1]$, $V:\R\times\Omega\to [0,1]$ are stationary ergodic processes satisfying the following assumptions, 
for some constant $\kappa > 0$:
\begin{itemize}
\item[(A)]  $\sqrt{a(\,\cdot\,,\omega)}:\R\to (0,1]$\ is $\kappa$--Lipschitz continuous for all $\omega\in\Omega$;\medskip
\item[(V)] $V(\,\cdot\,,\omega):\R\to [0,1]$\ is $\kappa$--Lipschitz
  continuous for all $\omega\in\Omega$.
\end{itemize}

As for the nonlinearity $G$, we assume that it belongs to the class $\Ham$ defined as follows. 
\begin{definition}\label{def:Ham}
	A function $G:\R\to\R$ is said to be in the class $\Ham$ if it satisfies the following conditions, for some constants $\alpha_0,\alpha_1>0$ and $\gamma>1$:
	\begin{itemize}
		\item[(G1)] $\alpha_0|p|^\gamma-1/\alpha_0\leqslant G(p)\leqslant\alpha_1(|p|^\gamma+1)$\ for all ${p\in\R}$;\medskip
		\item[(G2)] $|G(p)-G(q)|\leqslant\alpha_1\left(|p|+|q|+1\right)^{\gamma-1}|p-q|$\ for all $p,q\in\R$.
\end{itemize}
\end{definition}

Solutions, subsolutions and supersolutions of \eqref{eq:generalHJ}
will be always understood in the viscosity sense, see
\cite{users,barles_book,bardi}. Assumptions (A), (V), (G1) and (G2)
guarantee well-posedness in $\D{UC}(\ccyl)$ of the Cauchy problem for
the parabolic equation \eqref{eq:generalHJ}, as well as Lipschitz
estimates for the solutions under appropriate assumptions on the
initial condition, see Appendix \ref{app:PDE} for more
details.

\begin{remark}
  Note that a solution $u^\epsilon(t,x,\omega)$ of \eqref{eq:introHJ}
  can be obtained from a solution $u(t,x,\omega)$ of
  \eqref{eq:generalHJ} by a hyperbolic rescaling of time and space:
  $u^\epsilon(t,x,\omega)=\epsilon
  u\left(\frac{t}{\epsilon},\frac{x}{\epsilon},\omega\right)$. For
  this reason, homogenization of \eqref{eq:introHJ} can be
  considered as a law of large numbers type of result. The natural
  questions about fluctuations around the deterministic limit
  $\overline{u}(t,x)$ and large deviations remain wide open even in
  the case when $G(p)=\frac12p^2$ and $d=1$.\footnote{The known quantitative
  homogenization estimates we mentioned earlier
  (\cite{ACS14,AC15,AC18}) show that, under assumptions, the
  convergence rates are algebraic but little is known about the ``true
  size'' of fluctuations.}
\end{remark}

The purpose of this paper is to establish a homogenization result for
equation \eqref{eq:introHJ}. We succeeded in doing this under the
additional hypothesis that the pair $(a,V)$ satisfies the scaled hill
and valley condition, which we shall now define. We will say that
$(a,V)$ under $\P$ satisfies the {\em scaled hill} (respectively, {\em
  scaled valley}) {\em condition} if
\begin{itemize}
\item[(S)] for every $h\in (0,1)$ and $y>0$, there exists a set
  $\Omega(h,y)$ of probability 1 such that, for every
  $\omega\in \Omega(h,y)$, there exist $\ell_1<\ell_2$ (depending on
  $\omega$) such that
\begin{itemize}
\item{(a)} \quad $\displaystyle \int_{\ell_1}^{\ell_2} \frac{dx}{a(x,\omega)} = y$,
\end{itemize}
and 
\begin{itemize}
\item{(h)}\ {\em (hill)} \quad $V(\,\cdot\,,\omega)\geqslant h\quad\hbox{on $[\ell_1,\ell_2]$}$\medskip
\end{itemize}
(respectively,\smallskip
\begin{itemize}
	\item{(v)}\ {\em (valley)} \quad $V(\,\cdot\,,\omega)\leqslant h\quad\hbox{on $[\ell_1,\ell_2]$}$).\\
\end{itemize}
\end{itemize}

\begin{remark}\label{rem:vac}
  The scaled hill and valley condition can be adapted to the possibly
  degenerate case $a\ge 0$ by replacing $a$ in the integral of item
  (a) with $a\vee \delta$ for some $\delta=\delta(\omega)\in\ (0,1)$,
  see \cite[p.\,236]{DK22}. This more general formulation is
  equivalent to (S) when $a>0$, as is the case in the current
  paper. Moreover, when $a\equiv 0$, it reduces to merely requiring
  that $\essinf V(x,\omega) = 0$ and $\esssup V(x,\omega) = 1$.
\end{remark}

Our main result reads as follows.

\begin{theorem}\label{thm:genhom}
  Suppose $a$ and $V$ satisfy (A), (V) and (S), $\beta\ge 0$, and
  $G\in\Ham$. Then, the viscous HJ equation \eqref{eq:introHJ}
  homogenizes, i.e., there exists a continuous nonrandom function
  $\HV(G):\R\to\R$, {\em the effective Hamiltonian}, such that, with
  probability 1, for every uniformly continuous function $g$ on $\R$,
  the solutions $u^\epsilon$ of \eqref{eq:introHJ} satisfying
  $u^\epsilon(0,\,\cdot\,,\omega) = g$ converge locally uniformly on
  $[0, +\infty)\times \R$ as $\epsilon\to 0+$ to the unique solution
  $\ol{u}$ of the equation
  \[ \partial_t \ol{u} = \HV(G)(D\ol{u}), \quad (t,x)\in
    (0,+\infty)\times\R, \] satisfying $\ol{u}(0,\,\cdot\,) = g$.
\end{theorem}

\begin{remark}
  Our notation $\HV(G)$ shows only the dependence of the effective
  Hamiltonian on $G$. The inherent dependence of the effective
  Hamiltonian on the law of the pair of stochastic processes
  $(a(x,\,\cdot\,),V(x,\,\cdot\,))_{x\in\R}$ and on the constant
  $\beta$ is not reflected in the notation, since they will remain
  fixed throughout the paper. On the other hand, in our arguments
  leading to the proof of Theorem \ref{thm:genhom}, we will
  occasionally change $G$.
\end{remark}

\subsection{Proof strategy and outline of the paper}\label{sec:outline}

In this section, we will describe our strategy for proving the homogenization result stated in Theorem \ref{thm:genhom}. 
Let us denote by $u_\theta(\,\cdot\,,\,\cdot\,,\omega)$ the unique Lipschitz solution to \eqref{eq:generalHJ} with initial condition $u_\theta(0,x,\omega)=\theta x$ on $\R$, and introduce the following deterministic quantities, defined almost surely on $\Omega$:
\begin{eqnarray}\label{eq:infsup}
	\HV^L(G) (\theta):=\liminf_{t\to +\infty}\ \frac{u_\theta(t,0,\omega)}{t}\quad\text{and}\quad 
	\HV^U(G) (\theta):=\limsup_{t\to +\infty}\ \frac{u_\theta(t,0,\omega)}{t}.
\end{eqnarray}
In view of \cite[Lemma 4.1]{DK17}, for the purpose of proving homogenization, it suffices to show that $\HV^L(G)(\theta) = \HV^U(G)(\theta)$ for every $\theta\in\R$. In this case, the common value is denoted by $\HV(G)(\theta)$. The function $\HV(G):\R\to\R$ is called the effective Hamiltonian associated with $G$ and appears in the statement of Theorem \ref{thm:genhom}.

Note that the above observation readily implies homogenization when $\beta=0$. Indeed, in this instance, $u_\theta(t,x,\omega)=\theta x + tG(\theta)$ for all $(t,x)\in\ccyl$ and $\theta \in\R$, so $\HV(G) = G$. Hence, to prove our homogenization result, it suffices to consider $\beta>0$. 

In order to prove Theorem \ref{thm:genhom}, we will adapt the approach that was taken in \cite{KYZ20, DK22} and subsequently developed in \cite{Y21b}. It consists in showing the existence of suitable solutions for the (stationary) corrector equation associated to \eqref{eq:generalHJ}. Following \cite{Y21b}, the last problem reduces to considering the family of first-order ODEs  
\begin{equation}\label{eq:cellODE}
	a(x,\omega)f'(x)+G(f(x))+\beta V(x,\omega)=\lambda,\quad x\in\R,
\end{equation}
for every $\lambda \geqslant m_0 + \beta$ and $\omega\in\Omega$, where $m_0 := \inf\{G(p):\,p\in\R\}$. Note that, due to the assumptions (G1)--(G2), this infimum is attained and, in particular, it is finite.

We shall be interested in solutions to equation \eqref{eq:cellODE}
which satisfy the conditions stated in the next definition.  It will
be convenient to give a name to the subset of $\Ham$ for which such
solutions exists. Note that, unlike $\Ham$, this subset depends on the
choice of $a,\, V$ and $\beta$.

\begin{definition}\label{def:Hamone} 
  A function $G\in\Ham$ is said to be in the class
  $\Hamone(a,V,\beta)$ if, for every $\theta\in\mathbb{R}$, there
  exists a unique constant
  $\lambda = \lambda(\theta) \geqslant m_0 + \beta$ and stationary
  functions $\ul f,\ol f:\R\times\Omega\to\R$ (depending on $\theta$
  and not necessarily distinct) that satisfy the following conditions:
   \begin{itemize}
  \item[\em (i)] $\ul f(\,\cdot\,,\omega), \ol f(\,\cdot\,,\omega)$ are $\CC^1_b$ solutions of \eqref{eq:cellODE} on a set $\wh\Omega$ of probability 1;
  \item[\em (ii)] $\inf\{\ol f(x,\omega) - \ul f(x,\omega):\,x\in\R\} = 0$ on a
    set $\wh\Omega$ of probability $1$;
  \item[\em (iii)]
    $\mathbb{E}[\ul f(0,\omega)] \leqslant \theta \leqslant
    \mathbb{E}[\ol f(0,\omega)]$.
  \end{itemize}
\end{definition}

The class $\Hamone(a,V,\beta)$ is tailored to our proof of homogenization. In fact, the following holds. 

\begin{theorem}\label{thm:spehom}
  Suppose $a$ and $V$ satisfy (A) and (V), $\beta > 0$, and
  $G\in\Hamone(a,V,\beta)$. Then, the viscous HJ equation
  \eqref{eq:introHJ} homogenizes, and the effective Hamiltonian is
  given by
	\[ 
	\HV(G)(\theta) = \lambda(\theta)\quad\text{for all $\theta\in\R$}.
	\]
\end{theorem}

The proof of this theorem is given in Section \ref{sec:spehom}. We provide here an outline of the argument. 
As we noted at the beginning of this section, to prove homogenization, it
suffices to show that
$\HV^L(G)(\theta) = \HV^U(G)(\theta) = \lambda(\theta)$ for every
$\theta\in\R$. These equalities follow via the same argument as in \cite[Lemma
5.6]{DK22} if $\theta = \mathbb{E}[\ul f(0,\omega)]$ or
$\theta = \mathbb{E}[\ol f(0,\omega)]$. Hence, in the proof of Theorem
\ref{thm:spehom}, we will assume that
$\mathbb{E}[\ul f(0,\omega)] < \theta < \mathbb{E}[\ol
f(0,\omega)]$. For such $\theta$, every $\delta>0$ and $\P$-a.e.\
$\omega$, we carefully interpolate between $\ul f$ and $\ol f$ on a
suitable interval where $\ol f - \ul f \leqslant \delta$, and thereby
construct an $f_\delta$ such that
\[ a(x,\omega)f_\delta' + G(f_\delta) + \beta V(x,\omega) \geqslant \lambda - C_0\delta,\quad x\in\R, \]
for some $C_0 > 0$ (independent of $\delta$). We use $f_\delta$ to build a subsolution $v_\delta$ of \eqref{eq:generalHJ} such that
\[ \lim_{t\to+\infty}\frac{v_\delta(t,0,\omega)}{t} = \lambda(\theta) - C_0\delta\quad\text{and}\quad v_\delta(0,x,\omega) \leqslant u_\theta(0,x,\omega) = \theta x,\quad\forall x\in\R, \]
conclude that $\HV^L(G)(\theta) \geqslant \lambda(\theta) - C_0\delta$ by the comparison principle, and then let $\delta\to 0$. The inequality $\HV^U(G)(\theta) \leqslant \lambda(\theta)$ can be obtained similarly. This argument also shows that the effective Hamiltonian is constant on an interval of positive length that contains $\theta$ if and only if $\ul f,\ol f$ are almost surely distinct. 

The rest of the paper is devoted to extending the homogenization result from $\Hamone(a,V,\beta)$ to all of $\Ham$.  
To this end, we introduce the following new subclass of $\Ham$.

\begin{definition}\label{def:Hamtwo}
	A function $G\in\Ham$ is said to be in the class $\Hamzero$ if it satisfies the following additional conditions:
	\begin{itemize}
		\item[(G3)] $G(0) = m_0 = 0$ and $G(p) > 0$ for every $p\ne 0$;\smallskip
		\item[(G4)] $G$ has finitely many local extrema and it attains distinct values at these points.
	\end{itemize}
\end{definition}

The core of our paper consists in proving the next theorem. This is also the part where we need the scaled hill and valley condition (S).

\begin{theorem}\label{thm:induction}
	Suppose $a$ and $V$ satisfy (A), (V) and (S), and $\beta>0$. Then, $\Hamzero \subset \Hamone(a,V,\beta)$.
\end{theorem}

The broad structure of the proof of Theorem \ref{thm:induction} is parallel to the proof of homogenization in \cite{ATY_1d} in the inviscid setting, but the arguments involve new ideas and techniques due to the presence of the diffusive term. 
By a gluing procedure at the origin, which is provided in Section \ref{sub:glue1} (cf.\ \cite[Section 4.2]{ATY_1d}), the problem can be essentially reduced to the case when $G\in\Hamzero$ is strictly decreasing on $(-\infty,0]$. Next, suppose that $G$ has $2N+1$ local extrema on $[0,+\infty)$ for some $N\geqslant 1$, attained at $0 < p_1 < p_2 < \ldots < p_{2N}$. Let
\begin{equation}\label{eq:Mm}
\begin{aligned}
	&M = \max\{M_1,\ldots,M_N\} \quad\text{and}\quad m = \min\{m_1,\ldots,m_N\},\\
	&\text{where $M_i = G(p_{2i-1})$ and $m_j = G(p_{2j})$ are, respectively,}\\
	&\text{the local maxima and the local minima of $G$ on $(0,+\infty)$.}
\end{aligned}
\end{equation}
In Section \ref{sec:base}, we prove Theorem \ref{thm:induction} in the following two base cases:
\begin{itemize}
	\item [(i)] $G$ is strictly increasing on $[0,+\infty)$ (i.e., $N=0$ and $G$ quasiconvex on $\R$);\smallskip
	\item [(ii)] $G$ has multiple local extrema on $[0,+\infty)$, and $\beta > M - m$ (cf.\ \cite[Section 3]{ATY_1d}).
\end{itemize}

In the first base case, $G\in\Hamone(a,V,\beta)$ as a corollary of the results in \cite{Y21b} for quasiconvex Hamiltonians, and the effective Hamiltonian $\HV(G)$ is then also quasiconvex. In the second base case, we generalize the arguments in \cite{Y21b} to prove that $G\in\Hamone(a,V,\beta)$. Moreover, in that case, we show that $\HV(G)$ is quasiconvex even though $G$ is not. 

In Section \ref{sub:glue2}, we consider the case when $N\ge 1$ and $\beta \leqslant M - m$. First, we identify $\tilde G_1,\tilde G_2\in\Ham$ such that $G = \tilde G_1\wedge \tilde G_2$ and the number of local extrema of $\tilde G_1$ (respectively, $\tilde G_2$) is strictly less than $2N+1$. There are two subcases to analyze, and each one of them corresponds to a different choice of $\tilde G_1,\tilde G_2$, see Sections \ref{subsub:glue21}  and \ref{subsub:glue22}. Then, by further gluing arguments (cf.\ \cite[Sections 4.3--4.5]{ATY_1d}), we prove that, if $\tilde G_1,\tilde G_2\in\Hamone(a,V,\beta)$, then $G\in\Hamone(a,V,\beta)$ as well.

In Section \ref{sec:induction}, we give the proof of Theorem \ref{thm:induction}. It is based on a strong induction argument which puts together the information gathered in Sections \ref{sec:base} and \ref{sec:gluing}.

In Section \ref{sec:stable}, using a stability argument which is
provided in Appendix \ref{app:PDE}, we show that we can reduce the
proof of homogenization for $G\in\Ham$ to the case when
$G\in\Hamzero$.  The precise statement is the following.

\begin{theorem}\label{thm:main reduction}
	Suppose $a$ and $V$ satisfy (A) and (V), and $\beta>0$. Under these assumptions, if equation \eqref{eq:generalHJ} homogenizes for every $G\in\Hamzero$, then it homogenizes for every $G\in\Ham$. 
\end{theorem}

Clearly, our main homogenization result (i.e., Theorem \ref{thm:genhom}) is an immediate consequence of Theorems \ref{thm:spehom}, \ref{thm:induction} and \ref{thm:main reduction}. 

Some general results regarding the existence and properties of solutions of \eqref{eq:cellODE}, which crucially rely on the scaled hill and valley condition (S) and which are repeatedly used in Sections \ref{sec:base} and \ref{sec:gluing}, are presented in Section \ref{sec:lemmas}.

The paper ends with two appendices. In the first one, we have collected and proved several results on stationary sub and supersolutions of a class of equations that includes \eqref{eq:cellODE} as a particular case. The second one contains some PDE results needed for our proofs.

\section{Proof of Theorem \ref{thm:spehom}}\label{sec:spehom}

For every $\theta\in\R$, let
$\lambda = \lambda(\theta) \geqslant m_0 + \beta$ and $\ul f,\ol f$ be
as in Definition \ref{def:Hamone}. Recall the notation in
\eqref{eq:infsup}. By \cite[Lemma 4.1]{DK17}, it suffices to show that
	\[ \text{(i)}\quad\HV^L(G)(\theta)\geqslant \lambda \quad \text{and\quad (ii)} \quad \HV^U(G)(\theta)\leqslant \lambda. \]
	The proofs of (i) and (ii) are similar, so we will prove only the former.
	
        Recall from Definition \ref{def:Hamone} that
        $\mathbb{E}[\ul f(0,\omega)] \leqslant \theta \leqslant
        \mathbb{E}[\ol f(0,\omega)]$. If $\theta$ is an endpoint of
        this interval (which can be a singleton), then
        $\HV^L(G)(\theta) = \HV^U(G)(\theta) = \lambda$ by the same
        argument as in \cite[Lemma 5.6]{DK22}. Hence, we assume that
    \begin{equation}\label{eq:akatay}
    	\mathbb{E}[\ul f(0,\omega)] < \theta < \mathbb{E}[\ol f(0,\omega)].
    \end{equation}

    Let $\wh\Omega$ be as in Definition \ref{def:Hamone}. 	By Corollary \ref{cor:lelam} with $H(p,x,\omega) = G(p) + \beta V(x,\omega)$, there exists an $R>0$ such that
    \begin{equation}\label{eq:abur}
    	\max\{\|\ul f(\,\cdot\,,\omega)\|_\infty,\|\ol f(\,\cdot\,,\omega)\|_\infty\} \leqslant R
    \end{equation}
    on a set $\wh\Omega_1 \subset \wh\Omega$ of probability $1$. Fix a $y_0 \geqslant 1$. For every $\omega\in\wh\Omega_1$ and $\delta > 0$, by Lemma \ref{lem:charac} with $H(p,x,\omega) = G(p) + \beta V(x,\omega)$, there is an interval $[L_1,L_2]$ such that
    \[ \int_{L_1}^{L_2}\frac{dx}{a(x,\omega)} = y_0
    \quad\text{and}\quad
    \ol f(x,\omega) - \ul f(x,\omega) \leqslant \delta,\ \ \forall x\in[L_1,L_2].
    \]
    
    We construct a function $f_\delta(\,\cdot\,,\omega)\in\CC^1(\R)$ by interpolating $\ul f$ and $\ol f$ on $[L_1 - r,L_2 + r]$ for some $r>0$ (to be determined). Precisely, we set
	\begin{equation}\label{eq:zand}
		f_\delta(x,\omega)=\xi(x,\omega)\ul{f}(x,\omega)+(1-\xi(x,\omega))\ol{f}(x,\omega),
	\end{equation}
	where $\xi(\,\cdot\,,\omega)\in\CC^1(\R)$ is non-decreasing,
	\begin{equation}\label{eq:specs}
		\xi(x,\omega)\equiv 0\ \text{on}\ (-\infty, L_1-r] \quad \text{and} \quad \xi(x,\omega)\equiv 1\ \text{on}\ [L_2+r,+\infty).
	\end{equation}
	In other words, $\xi'(\,\cdot\,,\omega)$ is a continuous probability density function supported on $[L_1-r,L_2+r]$. Our goal is to find such a $\xi(\,\cdot\,,\omega)$ that also ensures the following: there exists a constant $C_0 > 0$ (independent of $\delta$) such that
	\begin{equation}\label{eq:ecza}
		a(x,\omega)f_\delta'(x,\omega) + G(f_\delta(x,\omega)) + \beta V(x,\omega) \geqslant \lambda - C_0\delta
	\end{equation}
	for every $x\in\R$ and $\P$-a.e.\ $\omega\in\wh\Omega_1$. This inequality will play a key role in the proof of (i).
	
	For every $\delta > 0$, we momentarily suppress $(x,\omega)$
        from the notation and observe that
	\begin{align}\label{eq:concur}
		af_\delta' + G(f_\delta) + \beta V &= \xi(a\ul{f}' + G(\ul{f}) + \beta V) + (1-\xi)(a\ol{f}' + G(\ol{f}) + \beta V)\nonumber\\
		&\quad + G(\xi\ul{f} + (1-\xi)\ol{f}) - \xi G(\ul{f}) - (1-\xi)G(\ol{f}) - a\xi'(\ol{f} - \ul{f})\nonumber\\
		&= \lambda + G(\xi\ul{f}+(1-\xi)\ol{f}) - \xi G(\ul{f}) - (1-\xi)G(\ol{f}) - a\xi'(\ol{f}-\ul{f}).
	\end{align}
	Recall \eqref{eq:abur}, denote the Lipschitz constant of $G$ on $[-R,R]$ by $C_R$, and let $K > 0$ be a constant. It is easy to see from \eqref{eq:specs} and \eqref{eq:concur} that, for every $\omega\in\wh\Omega_1$, if
	\begin{equation}\label{eq:ifal}
		(\ol{f} - \ul{f})(x,\omega) < 2\delta, \quad \forall x\in[L_1-r,L_2+r]
	\end{equation}
	and
	\begin{equation}\label{eq:husnu}
		a(x,\omega)\xi'(x,\omega) \leqslant K, \quad \forall x\in[L_1-r,L_2+r],
	\end{equation}
        then \eqref{eq:ecza} holds with $C_0 = 2(C_R + K)$ for every
        $x\in \R$.
    
    Suppose tentatively that we take $\xi(x,\omega)=\int_{-\infty}^x \zeta(s,\omega)\,ds$, where
    \begin{equation}\label{eq:zeta}
    	\zeta(x,\omega)=(y_0a(x,\omega))^{-1}\mathbbm{1}_{(L_1,L_2)}(x).
    \end{equation}
    Note that \eqref{eq:specs} and \eqref{eq:husnu} hold with $r = 0$ and $K = \frac1{y_0} \leqslant 1$. However, $\xi(\,\cdot\,,\omega)$ is not in $\CC^1(\R)$. We show below that a standard mollification argument resolves this issue.
    
	Since $a(\,\cdot\,,\omega) > 0$, the difference $(\ol{f}-\ul{f})(\,\cdot\,,\omega)$ is locally Lipschitz due to \eqref{eq:cellODE} and \eqref{eq:abur}. Hence, \eqref{eq:ifal} holds for some $r>0$. Take a standard sequence of even convolution kernels $\rho_n$ supported on $[-1/n,1/n]$. For all sufficiently large $n,\ n>1/r$, define $\zeta_n(x,\omega)=(\zeta(\,\cdot\,,\omega)*\rho_n)(x)$ with $\zeta(\,\cdot\,,\omega)$ as in \eqref{eq:zeta}, and let $\xi_n(x,\omega) = \int_{-\infty}^x\zeta_n(s,\omega)\,ds$. Note that $\xi_n(\,\cdot\,,\omega)$ satisfies \eqref{eq:specs}.
	Moreover, letting $\ul{a}(\omega) = \inf_{x\in[L_1-2r,L_2+2r]}a(x,\omega)>0$ and recalling that $y_0\geqslant 1$, the following inequalities hold:
	\begin{align*}
		a(x,\omega)\xi'_n(x,\omega)&\leqslant \frac{1}{y_0}\int_{-1/n}^{1/n}\frac{a(x,\omega)\rho_n(s)}{a(x+s,\omega)}\,ds\\
		&\leqslant 1+\int_{-1/n}^{1/n}\frac{(a(x,\omega)-a(x+s,\omega))\rho_n(s)}{a(x+s,\omega)}\,ds\leqslant 1+\frac{2\kappa}{\ul{a}(\omega) n}
	\end{align*}
	for all $x\in [L_1 - r,L_2 + r]$, where $\kappa$ is the Lipschitz constant of $\sqrt{a(\,\cdot\,,\omega)}$ from (A).\footnote{Note that $|a(x,\omega) - a(y,\omega)| \leqslant |\sqrt{a(x,\omega)} + \sqrt{a(y,\omega)}||\sqrt{a(x,\omega)} - \sqrt{a(y,\omega)}| \leqslant 2|\sqrt{a(x,\omega)} - \sqrt{a(y,\omega)}|$.}
	Choosing $n > \max\left\{\frac1r,\frac{2\kappa}{\ul{a}(\omega)}\right\}$, we deduce that $\xi_n(\,\cdot\,,\omega)$ satisfies \eqref{eq:husnu} with $K = 2$.
	Therefore, if we take $\xi(\,\cdot\,,\omega) = \xi_n(\,\cdot\,,\omega)$ in \eqref{eq:zand}, then \eqref{eq:ecza} holds with $C_0 = 2(C_R + 2)$ for every $x\in\R$ and $\omega\in\wh\Omega_1$.
	
	We are ready to prove (i). There is a set $\wh\Omega_2 \subset \wh\Omega_1$ of probability 1 on which the limits in \eqref{eq:infsup} hold. For every $\omega\in\wh\Omega_2$, let \[v_\delta(t,x,\omega)=(\lambda-C_0\delta)t+\int_0^xf_\delta(s,\omega)\,ds-M_\delta(\omega),\]
	where $M_\delta(\omega)$ will be chosen later to ensure that $v_\delta(0,x,\omega)\leqslant \theta x$ for all $x\in\R$. Note that $v_\delta(\,\cdot\,,\,\cdot\,,\omega)$ is a subsolution of \eqref{eq:generalHJ}. Indeed, 
	\begin{align*}
		&\partial_tv_\delta-a(x,\omega)\partial_{xx}^2v_\delta - G(\partial_x v_\delta)-\beta V(x,\omega)\\
		=\ &(\lambda-C_0\delta)-\left(a(x,\omega)f_\delta'(x,\omega)+G\left(f_\delta(x,\omega)\right)+\beta V(x,\omega)\right)\overset{\eqref{eq:ecza}}{\le} 0.
	\end{align*}
	By \eqref{eq:akatay} and the ergodic theorem, there is an $\wh\Omega_3 \subset \wh\Omega_2$ of probability $1$ such that, for every $\omega\in \wh\Omega_3$,
	\[\limsup_{x\to+\infty}\frac{v_\delta(0,x,\omega)}{x}=\lim_{x\to+\infty}\frac{1}{x}\int_{L_2+r}^x\ul{f}(s,\omega)\,ds=\E[\ul f(0,\omega)] < \theta\]
	and
	\[\liminf_{x\to-\infty}\frac{v_\delta(0,x,\omega)}{x}=\lim_{x\to-\infty}\frac{1}{|x|}\int_{x}^{L_1 - r}\ol{f}(s,\omega)\,ds=\E[\ol f(0,\omega)] > \theta.\]
	Therefore, for every $\omega\in\wh\Omega_3$, we can pick $M_\delta(\omega)$ large enough so that $v_\delta(0,x,\omega)\leqslant \theta x$ for all $x\in\R$.  By the comparison principle, $v_\delta(t,x,\omega)\leqslant u_\theta(t,x,\omega)$ on $\ccyl\times\wh\Omega_3$ and, hence,
	\[\HV^L(G)(\theta)=\liminf_{t\to+\infty}\frac{u_\theta(t,0,\omega)}{t}\geqslant\liminf_{t\to+\infty}\frac{v_\delta(t,0,\omega)}{t}=\lambda-C_0\delta\] 
	for every $\omega\in \wh\Omega_3$. Since $\delta > 0$ is arbitrary, this implies (i).

\section{Properties of solutions to equation \eqref{eq:cellODE}}\label{sec:lemmas}

In this section, we will highlight some significant properties of (deterministic) solutions of equation \eqref{eq:cellODE} and eventually prove the existence of a distinguished pair of stationary solutions of the same equation. These results are established for a function $G\in\Ham$ satisfying 
\[
G(\,\cdot\, + p_{\rm{min}}) - m_0 \in \Hamzero\ \ \text{for some $p_{\rm{min}}\in\R$,} 
\]
where $m_0 := \min\{G(p):\,p\in\R\} >-\infty$.\smallskip

For every $\lambda \geqslant m_0 + \beta$, let
\[ \mathcal{I}_\lambda(G) := \{ [p_1,p_2] \subset
  \R:\,\{G(p_1),G(p_2)\} = \{\lambda - \beta,\lambda\}\ \text{and}\
  \lambda - \beta < G(p) < \lambda\ \text{for all}\ p\in(p_1,p_2)
  \}. \] For every $[p_1,p_2] \in \mathcal{I}_\lambda(G)$ and
$\omega\in\Omega$, define
\[ \Sol(\lambda,[p_1,p_2],\omega) := \{f(\,\cdot\,,\omega)\in\CC^1(\R):\, \text{$f$ solves \eqref{eq:cellODE} and $p_1 \leqslant f(x,\omega) \leqslant p_2$ for all $x\in\R$}\}. 
\]
We start by noting the following.
\begin{lemma}\label{lem:kargo}
Let $\lambda \geqslant m_0 + \beta$ and $[p_1,p_2] \in \mathcal{I}_\lambda(G)$. For every 
$\omega\in\Omega$, the set $\Sol(\lambda,[p_1,p_2],\omega)$ is nonempty and compact in $\CC(\R)$. 
\end{lemma}

\begin{proof}
  Take two sequences $(\lambda_n)_{n\in\N}$ and $(\beta_n)_{n\in\N}$
  such that $\lambda_n\nearrow\lambda$ and
  $\lambda_n - \beta_n \searrow \lambda - \beta$ as $n\to
  +\infty$. Fix an arbitrary $\omega\in\Omega$ and apply Lemma \ref{lem:inbetw} with
  $H(p,x) = G(p) + \beta_n V(x,\omega) - \lambda_n$, $m(x) \equiv p_1$ and
  $M(x) \equiv p_2$ to deduce the existence of an $f_n\in\CC^1(\R)$
  such that
	\begin{equation}\label{eq:tunebeta}
		a(x,\omega)f_n'(x)+G(f_n(x))+\beta_n V(x,\omega) = \lambda_n,\quad x\in\R,
	\end{equation}
	and $p_1 < f_n(x) < p_2$ for all $x\in\R$. The sequence $(f_n)_{n\in\N}$ is equi-bounded. Moreover, it is locally equi-Lipschitz by \eqref{eq:tunebeta} since $a(\,\cdot\,,\omega) > 0$ on $\R$. Hence, by the Arzel\'a-Ascoli theorem, up to extracting a subsequence (not relabeled), it converges, locally uniformly on $\R$, to a continuous function $f$ that satisfies $p_1 \leqslant f(x) \leqslant p_2$ for all $x\in\R$. Again by \eqref{eq:tunebeta}, the derivatives $(f'_n)_{n\in\N}$ form a Cauchy sequence in $\CC(\R)$,  hence the functions $f_n$ actually converge to $f$ in the local $\CC^1$ topology and $f$ solves \eqref{eq:cellODE}, i.e., $f\in\Sol(\lambda,[p_1,p_2],\omega)\ne \emptyset$.  Via the same argument, one can show that any sequence 
$(f_n)_{n\in\N}$ in $\Sol(\lambda,[p_1,p_2],\omega)$ converges in $\CC(\R)$, up to a subsequence, to some 
$f\in \Sol(\lambda,[p_1,p_2],\omega)$. This proves the asserted compactness of $\Sol(\lambda,[p_1,p_2],\omega)$ .  
	\end{proof}

Our next results rely crucially on the scaled hill and valley condition (S). 

\begin{lemma}\label{lem:cerca1}
  Let $\lambda \geqslant m_0 + \beta$, $\delta > 0$, $y_0 > 0$ and
  $[p_1,p_2]\in\mathcal{I}_\lambda(G)$. For $\P$-a.e.\ $\omega$, the
  following holds:
	\begin{itemize}
		\item [(a)] there is an interval $[L_1,L_2]$ such that
		\[ \int_{L_1}^{L_2}\frac{dx}{a(x,\omega)} = y_0
		\quad\text{and}\quad
		p_1 \leqslant f(x) \leqslant p_1 + \delta\ \ \forall x\in[L_1,L_2]
		\]
		for every $f\in\Sol(\lambda,[p_1,p_2],\omega)$;\smallskip
		\item [(b)] 	there is an interval $[L_1,L_2]$ such that
		\[ \int_{L_1}^{L_2}\frac{dx}{a(x,\omega)} = y_0
		\quad\text{and}\quad
		p_2 - \delta \leqslant f(x) \leqslant p_2\ \ \forall x\in[L_1,L_2]
		\]
		for every $f\in\Sol(\lambda,[p_1,p_2],\omega)$.
	\end{itemize}
\end{lemma}

\begin{proof}
	(a) Suppose $G(p_1) = \lambda - \beta$ and $G(p_2) = \lambda$.\footnote{If $G(p_1) = \lambda$ and $G(p_2) = \lambda - \beta$, then we can work with
	\[ \check G(p) = G(-p),\quad \check a(x,\omega) = a(-x,\omega),\quad \check V(x,\omega) = V(-x,\omega)\quad\text{and}\quad \check f(x) = -f(-x). \]}
	Assume without loss of generality that $\delta\in(0,p_2-p_1)$ and that it is small enough to ensure that 
	\begin{equation}\label{eq:kucuk}
		G(p_1 + \delta) = \min\{G(p):\,p_1 + \delta \leqslant p \leqslant p_2\} < \lambda.
	\end{equation}
	Fix $h\in\left(1-\frac{G(p_1 + \delta) - (\lambda - \beta)}{2\beta},1\right)$ so that $\beta-G(p_1 + \delta) +\lambda <2\beta h$ or, equivalently,
	\begin{equation}\label{eq:hc}
		\lambda-\beta h - G(p_1 + \delta) < -\beta(1-h).
	\end{equation}
	Next, again without loss of generality, assume that
	\begin{equation}\label{eq:nohay}
		y_0 > \frac{p_2 - p_1}{\beta(1-h)}.
	\end{equation}
	By the scaled hill condition, there exists a set $\Omega(h,2y_0)$ of probability $1$ such that, for every $\omega\in\Omega(h,2y_0)$, there exist $\ell_1<\ell_2$ such that (S)(a) and (S)(h) hold with $y = 2y_0$. Let $L_1\in(\ell_1,\ell_2)$ be such that 
	\begin{equation}\label{eq:zeli}
		\int_{\ell_1}^{L_1}\frac{dx}{a(x,\omega)}=\int_{L_1}^{\ell_2}\frac{dx}{a(x,\omega)}=\frac{y}{2}=y_0.
	\end{equation}
	
	For every $\omega \in \Omega(h,2y_0)$ and $f\in\Sol(\lambda,[p_1,p_2],\omega)$, if $f(x) \geqslant p_1 + \delta$ for some $x\in[\ell_1,\ell_2]$, then
	\begin{equation}\label{eq:balta}
		f'(x) = \frac{\lambda-\beta V(x,\omega)-G(f(x))}{a(x,\omega)}\leqslant\frac{\lambda-\beta h-G(p_1 + \delta)}{a(x,\omega)} < -\frac{\beta (1-h)}{a(x,\omega)}
	\end{equation}
	by \eqref{eq:cellODE}, \eqref{eq:kucuk}--\eqref{eq:hc} and (S)(h).
	It follows that $p_1 \leqslant f(x_0) \leqslant p_1 + \delta$ for some $x_0\in[\ell_1,L_1]$. Indeed, otherwise
	\begin{align*}
		f(L_1) = f(\ell_1) + \int_{\ell_1}^{L_1} f'(x)\,dx & < f(\ell_1) - \beta(1-h)\int_{\ell_1}^{L_1} \frac{dx}{a(x,\omega)} \leqslant p_2 -\beta(1-h)y_0 < p_1
	\end{align*}
	by \eqref{eq:nohay}--\eqref{eq:balta}, contradicting the fact that $f\in\Sol(\lambda,[p_1,p_2],\omega)$. 
	Moreover, inequality \eqref{eq:balta} shows that $f'(x)<0$ at points $x\in [\ell_1,\ell_2]$ where $f(x)=p_1+\delta$. We infer 	
	that $p_1 \leqslant f(x) \leqslant p_1 + \delta$ for every $x\in[x_0,\ell_2] \supset [L_1,\ell_2]$, and the assertion follows by setting $L_2 := \ell_2$.\smallskip
	
	(b) Suppose $G(p_1) = \lambda - \beta$ and $G(p_2) = \lambda$ as in part (a). Assume without loss of generality that $\delta\in(0,p_2-p_1)$ and that it is small enough to ensure that 
	\begin{equation}\label{eq:buyuk}
		G(p_2 - \delta) = \max\{G(p):\,p_1 \leqslant p \leqslant p_2 - \delta\} > \lambda - \beta.
	\end{equation}
	Fix $h\in\left(0,\frac{\lambda - G(p_2 - \delta)}{2\beta}\right)$ so that
	\begin{equation}\label{eq:hc2}
		\lambda-\beta h - G(p_2 - \delta) > \beta h.
	\end{equation}
	Next, again without loss of generality, assume that
	\begin{equation}\label{eq:nohay2}
		y_0 > \frac{p_2 - p_1}{\beta h}.
	\end{equation}
	By the scaled valley condition, there exists a set $\Omega(h,2y_0)$ of probability $1$ such that, for every $\omega\in\Omega(h,2y_0)$, there exist $\ell_1<\ell_2$ such that (S)(a) and (S)(v) hold with $y = 2y_0$. Let $L_1\in(\ell_1,\ell_2)$ be such that 
	\begin{equation}\label{eq:zeli2}
		\int_{\ell_1}^{L_1}\frac{dx}{a(x,\omega)}=\int_{L_1}^{\ell_2}\frac{dx}{a(x,\omega)}=\frac{y}{2}=y_0.
	\end{equation}
	
	For every $\omega \in \Omega(h,2y_0)$ and $f\in\Sol(\lambda,[p_1,p_2],\omega)$, if $f(x) \leqslant p_2 - \delta$ for some $x\in[\ell_1,\ell_2]$, then
	\begin{equation}\label{eq:balta2}
		f'(x) = \frac{\lambda-\beta V(x,\omega)-G(f(x))}{a(x,\omega)} \geqslant \frac{\lambda-\beta h-G(p_2 - \delta)}{a(x,\omega)} > \frac{\beta h}{a(x,\omega)}
	\end{equation}
	by \eqref{eq:cellODE}, \eqref{eq:buyuk}--\eqref{eq:hc2} and (S)(v). It follows that $p_2-\delta \leqslant f(x_0) \leqslant p_2$ for some $x_0\in[\ell_1,L_1]$. Indeed, otherwise
	\begin{align*}
		f(L_1) = f(\ell_1) + \int_{\ell_1}^{L_1} f'(x)\,dx & > f(\ell_1) + \beta h\int_{\ell_1}^{L_1} \frac{dx}{a(x,\omega)}\geqslant p_1 + \beta hy_0 > p_2
	\end{align*}
	by \eqref{eq:nohay2}--\eqref{eq:balta2},  contradicting the fact that $f\in\Sol(\lambda,[p_1,p_2],\omega)$. 
	Moreover, inequality \eqref{eq:balta2} shows that $f'(x)>0$ at points $x\in [\ell_1,\ell_2]$ where $f(x)=p_2-\delta$. We infer 	
	that  $p_2 - \delta \leqslant f(x) \leqslant p_2$ for every $x\in[x_0,\ell_2] \supset [L_1,\ell_2]$, and the assertion follows by setting $L_2 := \ell_2$.
\end{proof}

\begin{lemma}\label{lem:cerca2}
Let $\lambda \geqslant m_0 + \beta$, $\delta > 0$, $y_0 > 0$ and $[p_1,p_2],[p_2,p_3]\in\mathcal{I}_\lambda(G)$.  
For $\P$-a.e.\ $\omega$, there is an interval $[L_1,L_2]$ such that
	\[ \int_{L_1}^{L_2}\frac{dx}{a(x,\omega)} = y_0
	\quad\text{and}\quad
	p_2 -\delta \leqslant f_1(x) \leqslant p_2 \leqslant f_2(x) \leqslant p_2 + \delta\ \ \forall x\in[L_1,L_2]
	\]
	for every $f_1\in\Sol(\lambda,[p_1,p_2],\omega)$ and $f_2\in\Sol(\lambda,[p_2,p_3],\omega)$.
\end{lemma}

\begin{proof}
	There are two cases: (i) $G(p_1) = G(p_3) = \lambda$ and $G(p_2) = \lambda - \beta$; (ii) $G(p_1) = G(p_3) = \lambda - \beta$ and $G(p_2) = \lambda$. We will treat the first case using the scaled hill condition. The second case can be treated similarly using the scaled valley condition (cf.\ parts (a) and (b) of Lemma \ref{lem:cerca1}).
	
	Assume without loss of generality that $\delta \in (0,(p_2 - p_1) \wedge (p_3 - p_2))$ and that it is small enough to ensure the following:
	\begin{equation}\label{eq:viyu}
	\begin{aligned}
		G(p_2 - \delta) &= \min\{G(p):\,p_1 \leqslant p \leqslant p_2 - \delta\} < \lambda;\\
		G(p_2 + \delta) &= \min\{G(p):\,p_2 + \delta \leqslant p \leqslant p_3\} < \lambda.
	\end{aligned}
	\end{equation}
	Fix $h\in\left(1 - \frac{\min\{G(p_2 - \delta),G(p_2 + \delta)\} - (\lambda - \beta)}{2\beta},1\right)$ so that
	\begin{equation}\label{eq:hczar}
		\lambda - \beta h - G(p_2 \pm \delta) < -\beta(1-h).
	\end{equation}
	Next, again without loss of generality, assume that
	\begin{equation}\label{eq:nohayzar}
		y_0 > \frac{(p_2 - p_1)\vee (p_3 - p_2)}{\beta(1-h)}.
	\end{equation}
	By the scaled hill condition, there exists a set $\Omega(h,3y_0)$ of probability $1$ such that, for every $\omega\in\Omega(h,3y_0)$, there exist $\ell_1<\ell_2$ such that (S)(a) and (S)(h) hold with $y = 3y_0$. Let $L_1,L_2\in(\ell_1,\ell_2)$ be such that 
	\begin{equation}\label{eq:zelizar}
		\int_{\ell_1}^{L_1}\frac{dx}{a(x,\omega)} = \int_{L_1}^{L_2}\frac{dx}{a(x,\omega)} = \int_{L_2}^{\ell_2}\frac{dx}{a(x,\omega)}=\frac{y}{3}=y_0.
	\end{equation}
	
	For every $\omega \in \Omega(h,3y_0)$ and $f_2\in\Sol(\lambda,[p_2,p_3],\omega)$, if $f_2(x) \geqslant p_2 + \delta$ for some $x\in[\ell_1,\ell_2]$, then
	\begin{equation}\label{eq:baltazar}
		f_2'(x) = \frac{\lambda - \beta V(x,\omega)-G(f_2(x))}{a(x,\omega)}\leqslant\frac{\lambda - \beta h - G(p_2 + \delta)}{a(x,\omega)} < -\frac{\beta (1-h)}{a(x,\omega)}
	\end{equation}
	by \eqref{eq:cellODE}, \eqref{eq:viyu}--\eqref{eq:hczar} and (S)(h).
	It follows that $p_2 \leqslant f_2(x_1) \leqslant p_2 + \delta$ for some $x_1\in[\ell_1,L_1]$. Indeed, otherwise
	\begin{align*}
		f_2(L_1) = f_2(\ell_1) + \int_{\ell_1}^{L_1} f_2'(x)\,dx & < f_2(\ell_1) - \beta(1-h)\int_{\ell_1}^{L_1} \frac{dx}{a(x,\omega)}\\ &\leqslant p_3 - \beta(1-h)y_0 < p_2
	\end{align*}
	by \eqref{eq:nohayzar}--\eqref{eq:baltazar}, contradicting the fact that $f_2\in\Sol(\lambda,[p_2,p_3],\omega)$. Moreover, it follows from \eqref{eq:baltazar} that $p_2 \leqslant f_2(x) \leqslant p_2 + \delta$ for every $x\in[x_1,\ell_2] \supset [L_1,L_2]$.
	
	Similarly, for every $\omega \in \Omega(h,3y_0)$ and $f_1\in\Sol(\lambda,[p_1,p_2],\omega)$, there exists an $x_2\in[L_2,\ell_2]$ such that $p_2 - \delta \leqslant f_1(x) \leqslant p_2$ for every $x\in[\ell_1,x_2] \supset [L_1,L_2]$. This can be shown by repeating the argument in the paragraph above with
	\[ \check G(p) = G(-p),\quad \check a(x,\omega) = a(-x,\omega),\quad \check V(x,\omega) = V(-x,\omega)\quad\text{and}\quad f_2(x) = -f_1(-x).\qedhere \]
\end{proof}

For every $\lambda \geqslant m_0 + \beta$, $[p_1,p_2] \in \mathcal{I}_\lambda(G)$, $x\in\R$ and $\omega\in\Omega$, let
\begin{equation}\label{eq:lazim}
\ul f(x,\omega):=\inf_{f\in\Sol(\lambda,[p_1,p_2],\omega)} f(x)
\quad\text{and}\quad
\ol f(x,\omega):=\sup_{f\in\Sol(\lambda,[p_1,p_2],\omega)} f(x).
\end{equation}
Note that $\ul f$ and $\ol f$ are not necessarily distinct functions.

\begin{lemma}\label{lem:measurable}
	The functions $\ul f, \ol f: \R\times\Omega\to\R$ are jointly measurable and stationary. Moreover, $\ul f(\,\cdot\,,\omega),\, \ol f(\,\cdot\,,\omega)\in \Sol(\lambda,[p_1,p_2],\omega)$ for every $\omega\in\Omega$.
\end{lemma}

\begin{proof}
The last assertion is a direct consequence of Lemmas \ref{lem:kargo} and \ref{appendix A lemma lattice}. The stationarity property of $\ul f$, $\ol f$ is a consequence of the fact that, for every fixed $z\in\R$ and $\omega\in\Omega$, $f\in \Sol(\lambda,[p_1,p_2],\omega)$ if and only if $f(\,\cdot\,+z)\in \Sol(\lambda,[p_1,p_2],\tau_z\omega)$ since $a(\,\cdot\,+z,\omega)=a(\,\cdot\,,\tau_z\omega)$ and $V(\,\cdot\, + z,\omega)=V(\,\cdot\,,\tau_z\omega)$. 

Let us prove that $\ul f:\R\times\Omega\to\R$ is measurable with respect to the product $\sigma$-algebra ${\mathcal B}\otimes\F$. This is equivalent to showing that $\Omega\ni \omega\mapsto \ul f(\,\cdot\,,\omega)\in\CC(\R)$ is a random variable from $(\Omega,\F)$ to the Polish space $\CC(\R)$ endowed with its Borel $\sigma$-algebra, see for instance  \cite[Proposition 2.1]{DS09}. 
Since the probability measure $\P$ is complete on $(\Omega,\F)$, it is enough to show that, for every fixed $\eps>0$, there exists a set $F\in \F$ with $\P(\Omega\setminus F)<\eps$ such that the restriction 
$\ul f$ to $F$ is a random variable from $F$ to $\CC(\R)$.  To this aim, we notice that the measure $\P$ is inner regular on $(\Omega,\F)$, see \cite[Theorem 1.3]{Bill99}, hence it is a Radon measure. By applying Lusin's Theorem \cite{LusinThm} to the random variables $a:\Omega\to\CC(\R)$ and $V:\Omega\to\CC(\R)$, we infer that there exists a compact set $F\subseteq \Omega$ with $\P(\Omega\setminus F)<\eps$ such that $a_{| F},V_{| F}:F\to\CC(\R)$ are continuous. We claim that $\ul f:\R\times F\to\R$ is lower semicontinuous. 
Indeed, let $(x_n,\omega_n)_{n\in\N}$ be a sequence converging to some $(x_0,\omega_0)$ in $\R\times F$. 
By the continuity of $a$ on $\R\times F$, we have that $\min_{J\times F} a>0$ for every compact interval $J\subset \R$. This implies that 
the functions  $\ul f(\,\cdot\,,\omega_n)$ are locally equi-Lipschitz on $\R$. Since they are also equi-bounded on $\R$ by the definition of $\Sol(\lambda,[p_1,p_2],\omega)$, by the Arzel\`a-Ascoli Theorem, we can extract a subsequence $\big(x_{n_k},\omega_{n_k}\big)_{k\in\N}$ such that $\liminf_n \ul f(x_n,\omega_n)=\lim_k \ul f(x_{n_k},\omega_{n_k})$ and 
$\ul f(\,\cdot\,,\omega_{n_k})$ converge to some $f$ in $\CC(\R)$. Since each function $\ul f(\,\cdot\,,\omega_{n})$ is a solution to \eqref{eq:cellODE} with $\omega:=\omega_n$ and $a(\,\cdot\,,\omega_n)\to a(\,\cdot\,,\omega_0)$, $V(\,\cdot\,,\omega_n)\to V(\,\cdot\,,\omega_0)$ in $\CC(\R)$, an argument analogous to the one used in the proof of Lemma  \ref{lem:kargo} shows that the functions $\ul f(\,\cdot\,,\omega_{n_k})$ actually converge to $f$ in the local $\CC^1$ topology. This readily implies that $f\in \Sol(\lambda,[p_1,p_2],\omega_0)$. By the definition of $\ul f$, we conclude that $\ul f(\,\cdot\,,\omega)\leqslant f$, in particular
\[
\ul f(x_0,\omega_0)
\leqslant
f(x_0)
=
\lim_k \ul f (x_{n_k},\omega_{n_k})
=
\liminf_n \ul f (x_n,\omega_{n}),
\]
proving the asserted lower semicontinuity property of $\underline f$. This implies that $\ul f_{|F}:F\to\CC(\R)$ is measurable (see, e.g., \cite[Proposition 2.1]{DS09}).
Via a similar argument, one can show that $\ol f:\R\times F\to\R$ is upper semicontinuous. 
\end{proof}

As far as stationary solutions to \eqref{eq:cellODE} are concerned, we have the following neat description.

\begin{lemma}\label{lem:clave}
	Take any $\lambda \geqslant m_0 + \beta$. Let $f(x,\omega)$ be a stationary process such that, for all $\omega\in\Omega$, $f(\,\cdot\,,\omega)\in\CC^1(\R)\cap \CC_b(\R)$ and it solves \eqref{eq:cellODE}. Then, there is an interval $[p_1,p_2] \in \mathcal{I}_\lambda(G)$ such that $\{f(x,\omega):\,x\in\R\} = (p_1,p_2)$ for $\P$-a.e.\ $\omega$.
	Consequently, $f(\,\cdot\,,\omega)\in \Sol(\lambda,[p_1,p_2],\omega)$ for $\P$-a.e.\ $\omega$.
\end{lemma}

\begin{proof}
First, take any $p_0\in\R$ such that $G(p_0) = \lambda - \beta$ (resp.\ $G(p_0) = \lambda$) and set 
$f_0(\,\cdot\,,\omega):= p_0$. 
 It is easy to check that the functions $f_1:= f_0$ and $f_2:= f$ (resp.\ $f_1:= f$ and $f_2:= f_0$) satisfy the assumptions of Lemma \ref{lem:order} with $H(p,x,\omega) = G(p) + \beta V(x,\omega)$.  Furthermore, $\P(f(x,\omega) = p_0\ \forall x\in\R) = 0$. We deduce from Lemma \ref{lem:order} that $\P(f(x,\omega) \neq p_0\ \forall x\in\R) = 1$, hence 
\begin{equation}\label{eq:duum}
	\P(\{ G(f(x,\omega)):\, x\in\R\} \cap \{\lambda - \beta,\lambda\} = \emptyset) = 1.
\end{equation}
Next, let us show that 
\begin{equation}\label{eq:teek}
	\P(\{ G(f(x,\omega)):\, x\in\R\} \supset (\lambda - \beta,\lambda) ) = 1.
	\end{equation}
To this aim, set $B:= \|f\|_\infty$. For every $\delta > 0$, pick $h\in(0,1)$ and $y > 0$ such that $2\beta(1-h) \leqslant \delta$ and $\beta (1-h)y > 2B$. By the scaled hill condition, for $\P$-a.e.\ $\omega$, there exist $\ell_1<\ell_2$ such that (S)(a) and (S)(h) hold. It follows that $G(f(x,\omega)) < \lambda - \beta + \delta$ for some $x\in[\ell_1,\ell_2]$. Indeed, otherwise, 
	\[ f'(x,\omega) = \frac{\lambda-\beta V(x,\omega)-G(f(x,\omega))}{a(x,\omega)} \leqslant \frac{\lambda-\beta h-(\lambda - \beta + \delta)}{a(x,\omega)} \leqslant \frac{\beta (1-h) - \delta}{a(x,\omega)} \leqslant -\frac{\beta (1-h)}{a(x,\omega)} \]
	for every $x\in[\ell_1,\ell_2]$ by \eqref{eq:cellODE}, yielding
	\[ f(\ell_2,\omega) = f(\ell_1,\omega) + \int_{\ell_1}^{\ell_2} f'(x,\omega)\,dx \leqslant f(\ell_1,\omega) - \beta(1-h)\int_{\ell_1}^{\ell_2} \frac{dx}{a(x,\omega)} \leqslant B -\beta(1-h)y < -B, \]
which is a contradiction. By an analogous argument involving the scaled valley condition, for every $\delta > 0$ and $\P$-a.e.\ $\omega$, there exists an $x\in\R$ such that $G(f(x,\omega)) > \lambda - \delta$. By continuity, we get  \eqref{eq:teek}. Finally, it follows from \eqref{eq:duum} and \eqref{eq:teek} that, there exists a set $\Omega_f$ of probability $1$ such that, for every $\omega\in\Omega_f$, there is an interval $[p_1,p_2] \in \mathcal{I}_\lambda(G)$ such that $\{f(x,\omega):\, x\in\R\} = (p_1,p_2)$. By ergodicity, we can take $\Omega_f$ such that the interval $[p_1,p_2]$ is the same for all $\omega\in\Omega_f$.
\end{proof}

\section{Base cases for Theorem \ref{thm:induction}}\label{sec:base}

\subsection{$G$ is quasiconvex}\label{sub:base1}

Take any $G\in\Hamzero$ such that
\[ G_1 := \left.G\right|_{(-\infty,0]}\ \text{is strictly decreasing\quad and}\quad G_2 := \left.G\right|_{[0,+\infty)}\ \text{is strictly increasing.} \]
Note that $G$ is quasiconvex, i.e., its sublevel sets $\{p\in\R:\, G(p) \leqslant \lambda\}$ are convex (intervals) for every $\lambda\geqslant 0$.

\begin{lemma}\label{lem:Y21b1}
	For every $\lambda \geqslant \beta$ and $\omega\in\Omega$, equation \eqref{eq:cellODE} has a unique solution $f_1^\lambda(\,\cdot\,,\omega)$ such that
	$f_1^\lambda(x,\omega)\in[G_1^{-1}(\lambda),G_1^{-1}(\lambda - \beta)]$ for all $x\in\R$.
	Similarly, it has a unique solution $f_2^\lambda(\,\cdot\,,\omega)$ such that $f_2^\lambda(x,\omega)\in[G_2^{-1}(\lambda - \beta),G_2^{-1}(\lambda)]$ for all $x\in\R$. Moreover, $f_1^\lambda$ and $f_2^\lambda$ are stationary.
\end{lemma}

\begin{proof}
	See \cite[Theorem 2.1]{Y21b}.
\end{proof}

Note that $[G_1^{-1}(\beta),0], [0,G_2^{-1}(\beta)] \in \mathcal{I}_\beta(G)$. By Lemma \ref{lem:cerca2}, for every $\delta > 0$, $y_0 > 0$ and $\P$-a.e.\ $\omega$, there is an interval $[L_1,L_2]$ such that
\[ \int_{L_1}^{L_2}\frac{dx}{a(x,\omega)} = y_0
\quad\text{and}\quad
-\delta \leqslant f_1^\beta(x,\omega) \leqslant 0 \leqslant f_2^\beta(x,\omega) \leqslant \delta\ \ \forall x\in[L_1,L_2].
\]

\begin{lemma}\label{lem:Y21b2}
	The quantities
	\[ \theta_1(\lambda) := \E[f_1^\lambda(0,\omega)] \quad\text{and}\quad \theta_2(\lambda) := \E[f_2^\lambda(0,\omega)] \]
	define two continuous bijections
	\[ \theta_1: [\beta,+\infty) \to (-\infty,\theta_1(\beta)] \quad\text{and}\quad \theta_2: [\beta,+\infty) \to [\theta_2(\beta),+\infty) \]
	which are decreasing and increasing, respectively.
\end{lemma}

	\begin{proof}
	See \cite[Theorem 2.2]{Y21b}.
\end{proof}

We recall the notation in Definition \ref{def:Hamone} and conclude that $G\in\Hamone(a,V,\beta)$ with the following choices:
\begin{itemize}
	\item if $\theta\in(-\infty,\theta_1(\beta))$, then $\lambda = \theta_1^{-1}(\theta) > \beta$ and $\ul f = \ol f = f_1^{\lambda}$;

	\item if $\theta\in[\theta_1(\beta),\theta_2(\beta)]$, then $\lambda = \beta$, $\ul f = f_1^\beta$ and $\ol f = f_2^{\beta}$;
	
	\item if $\theta\in(\theta_2(\beta),+\infty)$, then $\lambda = \theta_2^{-1}(\theta) > \beta$ and $\ul f = \ol f = f_2^{\lambda}$.
\end{itemize}

\subsection{$\beta$ is large}\label{sub:base2}

Take any $G\in\Hamzero$ such that
\begin{align*}
	&G_1 := \left.G\right|_{(-\infty,0]} \text{is strictly decreasing, and}\\
	&G_2 := \left.G\right|_{[0,+\infty)}\ \text{has multiple local extrema, attained at}\ 0 < p_1 < p_2 < \ldots < p_{2N}.
\end{align*}
Let $M = \max\{M_1,\ldots,M_N\}$ and $m = \min\{m_1,\ldots,m_N\}$ as
in \eqref{eq:Mm}.  Assume that
\[ \beta > M - m. \]

Define two right inverses $\ul G_2^{-1}, \ol G_2^{-1}: [0,+\infty) \to [0,+\infty)$ of $G_2$ by setting
\[\ul G_2^{-1}(\lambda) = \min\{p\geqslant 0:\, G(p) = \lambda\}\quad \text{and}\quad \ol G_2^{-1}(\lambda) = \max\{p\geqslant 0:\, G(p) = \lambda\} \]
for every $\lambda \geqslant 0$. Note that $\ul G_2^{-1}$ is strictly increasing, left-continuous everywhere, with discontinuities at some (but not necessarily all) $M_i$'s. Similarly, $\ol G_2^{-1}$ is strictly increasing, right-continuous everywhere, with discontinuities at some $m_j$'s.

\begin{lemma}
	$\ol G_2^{-1}(\lambda - \beta) < \ul G_2^{-1}(\lambda)$ and $[\ol G_2^{-1}(\lambda - \beta), \ul G_2^{-1}(\lambda)] \in \mathcal{I}_\lambda(G)$ for every $\lambda \geqslant \beta$.
\end{lemma}

\begin{proof}
	To justify the first claim, since $\beta > M - m$, it suffices to note the following:
	\begin{itemize}
		\item if $\lambda < m + \beta$, then $\ol G_2^{-1}(\lambda - \beta) = \ul G_2^{-1}(\lambda - \beta) < \ul G_2^{-1}(\lambda)$;
		\item if $\lambda > M$, then $\ol G_2^{-1}(\lambda - \beta) < \ol G_2^{-1}(\lambda) = \ul G_2^{-1}(\lambda)$.
	\end{itemize}
	Since $G(0) = 0$ and $\displaystyle{\lim_{p\to +\infty}G(p) = +\infty}$, the following hold, justifying the second claim:
	\begin{itemize}
		\item $G(p) > \lambda - \beta$ for every $p > \ol G_2^{-1}(\lambda - \beta)$;
		\item $G(p) < \lambda$ for every $p < \ul G_2^{-1}(\lambda)$. \qedhere
	\end{itemize}
\end{proof}

We apply Lemma \ref{lem:kargo} and deduce that $\Sol(\lambda,[\ol G_2^{-1}(\lambda - \beta), \ul G_2^{-1}(\lambda)],\omega) \ne \emptyset$ for every $\omega\in\Omega$. Let
\begin{align*}
&\ul f_2^\lambda(x,\omega):=\inf\{f(x):\,f\in\Sol(\lambda,[\ol G_2^{-1}(\lambda - \beta), \ul G_2^{-1}(\lambda)],\omega)\} \quad\text{and}\\
&\ol f_2^\lambda(x,\omega):=\sup\{f(x):\,f\in\Sol(\lambda,[\ol G_2^{-1}(\lambda - \beta), \ul G_2^{-1}(\lambda)],\omega)\}
\end{align*}
as in \eqref{eq:lazim}. By Lemma \ref{lem:measurable}, these functions are jointly measurable and stationary. Moreover, $\ul f_2^\lambda(\,\cdot\,,\omega),\, \ol f_2^\lambda(\,\cdot\,,\omega)\in \Sol(\lambda,[\ol G_2^{-1}(\lambda - \beta), \ul G_2^{-1}(\lambda)],\omega)$ for every $\omega\in\Omega$.

\begin{lemma}\label{lem:mon}
	If $\beta\leqslant\lambda_1<\lambda_2$, then $\P((\ol f_2^{\lambda_1} - \ul f_2^{\lambda_2})(x,\omega) < 0\ \ \forall x\in\R) = 1$.
\end{lemma}

\begin{proof}
	Pick $0 < \delta < \ol G_2^{-1}(\lambda_2 - \beta) - \ol G_2^{-1}(\lambda_1 - \beta)$. By Lemma \ref{lem:cerca1}(a), for $\P$-a.e.\ $\omega$, there is an $x\in\R$ such that
	\[ \ol f_2^{\lambda_1}(x,\omega) \leqslant \ol G_2^{-1}(\lambda_1 - \beta) + \delta < \ol G_2^{-1}(\lambda_2 - \beta) \leqslant \ul f_2^{\lambda_2}(x,\omega). \footnote{Alternatively, 
	pick $0 < \delta < \ul G_2^{-1}(\lambda_2) - \ul G_2^{-1}(\lambda_1)$. By Lemma \ref{lem:cerca1}(b), for $\P$-a.e.\ $\omega$, there is an $x\in\R$ such that
	$ \ol f_2^{\lambda_1}(x,\omega) \leqslant \ul G_2^{-1}(\lambda_1) < \ul G_2^{-1}(\lambda_2) - \delta \leqslant \ul f_2^{\lambda_2}(x,\omega)$.}\]
	Therefore, 
	$\P((\ol f_2^{\lambda_1} - \ul f_2^{\lambda_2})(x,\omega) < 0\ \text{for some}\ x\in\R) = 1$. The desired result now follows from Lemma \ref{lem:order} with $H(p,x,\omega) = G(p) + \beta V(x,\omega)$.
\end{proof}

\begin{lemma}\label{lem:con}
	For every $\lambda \geqslant \beta$ and $\lambda_n \geqslant \beta$, $n\in\N$, if $\lambda_n\nearrow\lambda$, then, for $\P$-a.e.\ $\omega$, as $n\to+\infty$, the sequence $\ul f_2^{\lambda_n}(\,\cdot\,,\omega)$ converges locally uniformly on $\R$ to $\ul f_2^\lambda(\,\cdot\,,\omega)$.
	Similarly, if $\lambda_n\searrow\lambda$, then, for $\P$-a.e.\ $\omega$, as $n\to+\infty$, the sequence $\ol f_2^{\lambda_n}(\,\cdot\,,\omega)$ converges locally uniformly on $\R$ to $\ol f_2^\lambda(\,\cdot\,,\omega)$.
\end{lemma}

\begin{proof}
	If $\lambda_n\nearrow\lambda$, then, by an Arzel\'a-Ascoli argument (as in the proof of Lemma \ref{lem:kargo}), for every $\omega\in\Omega$, the sequence $\ul f_2^{\lambda_n}(\,\cdot\,,\omega)$ has a subsequence that converges locally uniformly on $\R$ to some $f(\,\cdot\,,\omega)$ that solves \eqref{eq:cellODE} and $\ol G_2^{-1}((\lambda - \beta) -) \leqslant f(x,\omega) \leqslant \ul G_2^{-1}(\lambda-)$ for all $x\in\R$. Note that
	\[ \ol G_2^{-1}((\lambda - \beta) -) \leqslant \ol G_2^{-1}(\lambda - \beta) < \ul G_2^{-1}(\lambda) = \ul G_2^{-1}(\lambda-), \]
	and $[\ol G_2^{-1}(\lambda - \beta),\ul G_2^{-1}(\lambda)]$ is the unique element of $\mathcal{I}_\lambda(G)$ that is a subset of $[\ol G_2^{-1}((\lambda - \beta)-),\ul G_2^{-1}(\lambda)]$. We apply Lemma \ref{lem:clave} and deduce that $f(\,\cdot\,,\omega) \in \Sol(\lambda,[\ol G_2^{-1}(\lambda - \beta), \ul G_2^{-1}(\lambda)],\omega)$ for $\P$-a.e.\ $\omega$.
	
	By Lemma \ref{lem:mon}, for $\P$-a.e.\ $\omega$, $\ul f_2^{\lambda_n}(\,\cdot\,,\omega) < \ul f_2^{\lambda_{n+1}}(\,\cdot\,,\omega) < \ul f_2^\lambda(\,\cdot\,,\omega)$ for all $n\in\N$.
	Therefore, there is no need to pass to a subsequence, and $f(\,\cdot\,,\omega)\leqslant \ul f_2^\lambda(\,\cdot\,,\omega)$.
	By the minimality of $\ul f_2^\lambda(\,\cdot\,,\omega)$ in $\Sol(\lambda,[\ol G_2^{-1}(\lambda - \beta), \ul G_2^{-1}(\lambda)],\omega)$, we deduce that $f(\,\cdot\,,\omega) = \ul f_2^\lambda(\,\cdot\,,\omega)$.
	This concludes the proof of the first statement. The second statement is proved similarly.
\end{proof}

Recall that $G_1 := \left.G\right|_{(-\infty,0]}$ is strictly decreasing. For every $\lambda\geqslant \beta$, let $f_1^\lambda$ be as in Lemma \ref{lem:Y21b1}. Note that $[G_1^{-1}(\beta),0] \in \mathcal{I}_\beta(G)$ and $f_1^\beta(\,\cdot\,,\omega) \in \Sol(\beta,[G_1^{-1}(\beta),0],\omega)$ for every $\omega\in\Omega$. By Lemma \ref{lem:cerca2}, for every $\delta > 0$, $y_0 > 0$ and $\P$-a.e.\ $\omega$, there is an interval $[L_1,L_2]$ such that
\[ \int_{L_1}^{L_2}\frac{dx}{a(x,\omega)} = y_0
\quad\text{and}\quad
-\delta \leqslant f_1^\beta(x,\omega) \leqslant 0 \leqslant \ol f_2^\beta(x,\omega) \leqslant \delta\ \ \forall x\in[L_1,L_2].
\]
Similarly, by Lemma \ref{lem:cerca1}, for every $\lambda > \beta$, $\delta > 0$, $y_0 > 0$ and $\P$-a.e.\ $\omega$, there is an interval $[L_1,L_2]$ such that
\[ \int_{L_1}^{L_2}\frac{dx}{a(x,\omega)} = y_0
\quad\text{and}\quad
\ol f_2^\lambda(x,\omega) - \ul f_2^\lambda(x,\omega) \leqslant \delta\ \ \forall x\in[L_1,L_2].
\]

For every $\lambda \geqslant \beta$, define $\theta_1(\lambda) := \E[f_1^\lambda(0,\omega)]$ as in Lemma \ref{lem:Y21b2}. We introduce two generalizations of the quantity $\theta_2(\lambda)$ from the same lemma.

\begin{lemma}
	The quantities
	\[ \ul \theta_2(\lambda) := \E[\ul f_2^\lambda(0,\omega)] \quad\text{and}\quad \ol \theta_2(\lambda) := \E[\ol f_2^\lambda(0,\omega)] \]
	define two strictly increasing functions
	\[ \ul \theta_2: [\beta,+\infty) \to [\ul \theta_2(\beta),+\infty) \quad\text{and}\quad \ol \theta_2: [\beta,+\infty) \to [\ol \theta_2(\beta),+\infty) \]
	which are left-continuous and right-continuous, respectively. Moreover,
	\[ \bigcup_{\lambda\in(\beta,+\infty)}[\ul \theta_2(\lambda),\ol \theta_2(\lambda)] = (\ol \theta_2(\beta),+\infty). \]
\end{lemma}

\begin{proof}
	These statements follow immediately from Lemmas \ref{lem:mon}--\ref{lem:con} and the bounded convergence theorem.
\end{proof}

We recall the notation in Definition \ref{def:Hamone} and conclude that $G\in\Hamone(a,V,\beta)$ with the following choices:
\begin{itemize}
	\item if $\theta\in(-\infty,\theta_1(\beta))$, then $\lambda = \theta_1^{-1}(\theta) > \beta$ and $\ul f = \ol f = f_1^{\lambda}$;
	
	\item if $\theta\in[\theta_1(\beta),\ol \theta_2(\beta)]$, then $\lambda = \beta$, $\ul f = f_1^\beta$ and $\ol f = \ol f_2^{\beta}$;
	
	\item if $\theta\in(\ol \theta_2(\beta),+\infty)$, then $\theta\in[\ul \theta_2(\lambda),\ol \theta_2(\lambda)]$ for some $\lambda > \beta$, $\ul f = \ul f_2^{\lambda}$ and $\ol f = \ol f_2^{\lambda}$.\\
\end{itemize}

\section{Reduction steps for Theorem \ref{thm:induction}}\label{sec:gluing}

\subsection{Gluing at the origin}\label{sub:glue1}

\begin{figure}[h!]
  \begin{tikzpicture}
  \begin{axis}[scale only axis=true,
        width=0.6\textwidth,
        height=0.25\textwidth,
        axis x line=middle, xlabel={$p$}, axis y line=middle,
    ylabel={$\ $}, tick align=outside, samples=100,
    xtick={0},
    xticklabels={0}, ytick={0}, xmin=-3.8, xmax=3.5, ymin=0.01, ymax=3.3, very thick]
    \addplot[smooth, no marks, blue, solid] coordinates {
      (-3.5,3.0)
      (-3,1.3)
      (-2.3,2.5)
      (-1.9,1.7)
      (-1.5,2)
      (-0.5,0.35)
      (0,0)
      (0.5,0.3)
      (1.5,1.5)
      (2.3,1)
      (3,3)
    };
    \addlegendentry{$G$}
    \addplot+[smooth, no marks, black, dotted] coordinates {
      (-3.5,3.0)
      (-3,1.3)
      (-2.3,2.5)
      (-1.9,1.7)
      (-1.5,2)
      (-0.5,0.35)
      (0,0)
      (0.5,0.5)
      (1,2)
    };
    \addlegendentry{$\tilde{G}_1$}
    
    \addplot+[smooth, no marks, red, dashed] coordinates {
      (-1,2)
      (-0.5,0.5)
      (0,0)
      (0.5,0.3)
      (1.5,1.5)
      (2.3,1)
      (3,3)
    };
    \addlegendentry{$\tilde{G}_2$}
\end{axis}
\end{tikzpicture}
\caption{$G=\tilde{G}_1\wedge\tilde{G}_2$.}
\label{fig1}
\end{figure}
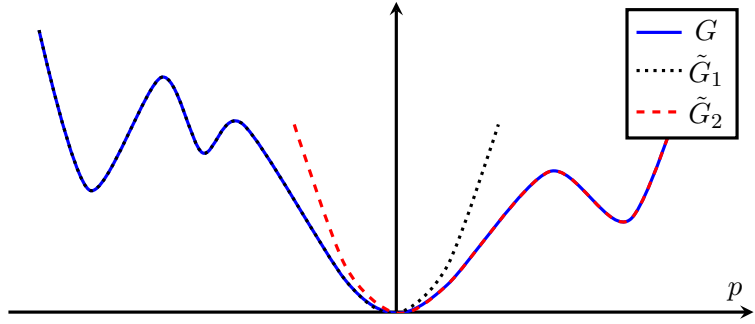

Take any $G\in\Hamzero$. Define
\[
\tilde G_1(p) = \begin{cases} G(p) &\text{if}\ p\leqslant 0,\\ I(p) &\text{if}\ p\geqslant 0,\end{cases}
\]
and
\[
\tilde G_2(p) = \begin{cases} I(p) &\text{if}\ p\leqslant 0,\\G(p) &\text{if}\ p\geqslant 0,\end{cases}
\]
where $I:\R\to[0,+\infty)$ is chosen such that
\[ I_1 := \left.I\right|_{(-\infty,0]}\ \text{is strictly
    decreasing\quad and}\quad I_2 := \left.I\right|_{[0,+\infty)}\
  \text{is strictly increasing,} \] $\tilde G_1, \tilde G_2\in\Ham_0$
and $G = \tilde G_1\wedge \tilde G_2$ (see Figure~\ref{fig1}). 

We will prove that, if $\tilde G_1,\tilde G_2\in\Hamone(a,V,\beta)$, then
$G\in\Hamone(a,V,\beta)$. Let
\[
\ol G_-^{-1}(\beta) = \max\{p\leqslant 0:\, G(p) = \beta \} \quad\text{and}\quad \ul G_+^{-1}(\beta) = \min\{p\geqslant 0:\, G(p) = \beta \}.
\]
Note that $[\ol G_-^{-1}(\beta),0], [0,\ul G_+^{-1}(\beta)] \in \mathcal{I}_\beta(G)$. We apply Lemma \ref{lem:kargo} and deduce that, for every $\omega\in\Omega$, $\Sol(\beta,[\ol G_-^{-1}(\beta),0],\omega) \ne \emptyset$ and $\Sol(\beta,[0,\ul G_+^{-1}(\beta)],\omega) \ne \emptyset$. Let
\begin{align*}
	&\ul f_-^\beta(x,\omega):=\inf\{f(x):\,f\in\Sol(\beta,[\ol G_-^{-1}(\beta),0],\omega)\} \quad\text{and}\\
	&\ol f_+^\beta(x,\omega):=\sup\{f(x):\,f\in\Sol(\beta,[0,\ul G_+^{-1}(\beta)],\omega)\}
\end{align*}
as in \eqref{eq:lazim}. By Lemma \ref{lem:measurable}, these functions are jointly measurable and stationary. Moreover, $\ul f_-^\beta(\,\cdot\,,\omega)\in\Sol(\beta,[\ol G_-^{-1}(\beta),0],\omega)$ and $\ol f_+^\beta(\,\cdot\,,\omega)\in\Sol(\beta,[0,\ul G_+^{-1}(\beta)],\omega)$ for every $\omega\in\Omega$. By Lemma \ref{lem:cerca2}, for every $\delta > 0$, $y_0 > 0$ and $\P$-a.e.\ $\omega$, there is an interval $[L_1,L_2]$ such that
\[ \int_{L_1}^{L_2}\frac{dx}{a(x,\omega)} = y_0
\quad\text{and}\quad
-\delta \leqslant \ul f_-^\beta(x,\omega) \leqslant 0 \leqslant \ol f_+^\beta(x,\omega) \leqslant \delta\ \ \forall x\in[L_1,L_2].
\]

Let
\begin{equation}
  \label{flat0}
  \ul\theta_-(\beta) = \mathbb{E}[\ul f_-^\beta(0,\omega)] < 0 \quad\text{and}\quad \ol\theta_+(\beta) = \mathbb{E}[\ol f_+^\beta(0,\omega)] > 0.
\end{equation}
Fix any $\theta > \ol\theta_+(\beta)$. Recall that $\tilde G_2\in\Hamzero\cap\Hamone(a,V,\beta)$ by our assumptions. Take $\lambda = \lambda(\theta)\geqslant\beta$ and $\ul f,\ol f$ as in Definition \ref{def:Hamone} for $\tilde G_2$. Since $G = \tilde G_2$ on $[0,+\infty)$, it follows from the lemma below that $\ul f,\ol f$ are stationary solutions  of \eqref{eq:cellODE} with $G$.

\begin{lemma}\label{lem:bound1}
	With the notation in the paragraph above, $\lambda > \beta$ and $\ol f(\,\cdot\,,\omega) \geqslant \ul f(\,\cdot\,,\omega) > \ol f_+^\beta(\,\cdot\,,\omega) \geqslant 0$ on $\mathbb{R}$ for $\P$-a.e.\ $\omega$.
\end{lemma}

\begin{proof}
	Suppose (for the purpose of reaching a contradiction) that $\lambda = \beta$. By Lemma \ref{lem:clave} (with $\tilde G_2$), for $\P$-a.e.\ $\omega$, $\ol f(\,\cdot\,,\omega)$ stays in an interval in $\mathcal{I}_\beta(\tilde G_2)$. There are only two such intervals: $[I_1^{-1}(\beta),0]$ and $[0,\ul G_+^{-1}(\beta)]$. It follows that, for $\P$-a.e.\ $\omega$,
	\begin{equation}\label{eq:roll}
		\ol f(\,\cdot\,,\omega) \leqslant \ol f_+^\beta(\,\cdot\,,\omega)\ \text{on}\ \R.
	\end{equation}
	Indeed, \eqref{eq:roll} is trivially true if $\ol f(\,\cdot\,,\omega)$ stays in $[I_1^{-1}(\beta),0]$. If it stays in $[0,\ul G_+^{-1}(\beta)]$, then \eqref{eq:roll} is true since $G = \tilde G_2$ on $[0,+\infty)$ and $\ol f_+^\beta(\,\cdot\,,\omega)$ is maximal in $\Sol(\beta,[0,\ul G_+^{-1}(\beta)],\omega)$. This implies that $\theta \leqslant \mathbb{E}[\ol f(0,\omega)] \leqslant \mathbb{E}[\ol f_+^\beta(0,\omega)] = \ol\theta_+(\beta) < \theta$, which is a contradiction. Therefore, $\lambda > \beta$.
	
	By Corollary \ref{cor:ordqu} with $H(p,x,\omega) = G(p) + \beta V(x,\omega)$, for $\P$-a.e.\ $\omega$, the distance between $\ol f_+^\beta(\,\cdot\,,\omega)$ and $\ol f(\,\cdot\,,\omega)$ (resp.\ $\ul f(\,\cdot\,,\omega)$) is greater than $\delta(\lambda - \beta)$ for some $\delta>0$. Since $\mathbb{E}[\ol f_+^\beta(0,\omega)] = \ol\theta_+(\beta) < \theta \leqslant \mathbb{E}[\ol f(0,\omega)]$, we see that $\ol f_+^\beta(\,\cdot\,,\omega) + \delta(\lambda - \beta) \leqslant \ol f(\,\cdot\,,\omega)$ on $\mathbb{R}$ for $\P$-a.e.\ $\omega$. Recalling from Definition \ref{def:Hamone} that $\inf\{\ol f(x,\omega) - \ul f(x,\omega):\,x\in\R\} = 0$ for $\P$-a.e.\ $\omega$, we conclude that $\ol f_+^\beta(\,\cdot\,,\omega) + \delta(\lambda - \beta) \leqslant \ul f(\,\cdot\,,\omega)$ on $\mathbb{R}$ for $\P$-a.e.\ $\omega$.
\end{proof}

Similarly, fix any $\theta < \ul\theta_-(\beta)$. Recall that $\tilde G_1\in\Hamzero\cap\Hamone(a,V,\beta)$ by our assumptions. Take $\lambda = \lambda(\theta)\geqslant\beta$ and $\ul f,\ol f$ as in Definition \ref{def:Hamone} for $\tilde G_1$. Since $G = \tilde G_1$ on $(-\infty,0]$, it follows from the lemma below that $\ul f,\ol f$ are stationary solutions  to \eqref{eq:cellODE} with $G$.

\begin{lemma}\label{lem:bound2}
	With the notation in the paragraph above, $\lambda > \beta$ and $\ul f(\,\cdot\,,\omega) \leqslant \ol f(\,\cdot\,,\omega) < \ul f_-^\beta(\,\cdot\,,\omega) \leqslant 0$ on $\mathbb{R}$ for $\P$-a.e.\ $\omega$.
\end{lemma}

\begin{proof}
	Analogous to the proof of Lemma \ref{lem:bound1}.
\end{proof}

We recall the notation in Definition \ref{def:Hamone} and conclude that $G\in\Hamone(a,V,\beta)$ with the following choices:
\begin{itemize}
	\item if $\theta\in(-\infty,\ul\theta_-(\beta))$, then $\lambda > \beta$ and $\ul f,\ol f$ are as in Definition \ref{def:Hamone} for $\tilde G_1$;
	
	\item if $\theta\in[\ul\theta_-(\beta),\ol\theta_+(\beta)]$, then $\lambda = \beta$, $\ul f = \ul f_-^\beta$ and $\ol f = \ol f_+^\beta$; moreover, $\ul\theta_-(\beta)<\ol\theta_+(\beta)$ by \eqref{flat0};
	
	\item if $\theta\in(\ol\theta_+(\beta),+\infty)$, then $\lambda > \beta$ and $\ul f,\ol f$ are as in Definition \ref{def:Hamone} for $\tilde G_2$.
\end{itemize}

\subsection{Gluing when $\beta$ is small}\label{sub:glue2}

Take any $G\in\Hamzero$ such that
\begin{align*}
	&\text{$G$ is strictly decreasing on $(-\infty,0]$, and}\\
	&\text{$G$ has multiple local extrema on $[0,+\infty)$, attained at $0 < p_1 < p_2 < \ldots < p_{2N}$.}
\end{align*}
Let $M = \max\{M_1,\ldots,M_N\}$ and $m = \min\{m_1,\ldots,m_N\}$ as
in \eqref{eq:Mm}.  By Definition \ref{def:Hamtwo},
\begin{equation}\label{eq:ij}
	\text{there is a unique $i\in\{1,\ldots,N\}$ (resp.\ $j\in\{1,\ldots,N\}$) such that $M = M_i$ (resp.\ $m = m_j$).}
\end{equation}
Assume that
\[ \beta \leqslant M - m. \]
There are two subcases to consider:
\begin{itemize}
	\item [(i)] $\beta > M - \min\{m_i,\ldots,m_N\}$;
	\item [(ii)] $\beta \leqslant M - \min\{m_i,\ldots,m_N\}$.
\end{itemize}

\subsubsection{$\beta > M - \min\{m_i,\ldots,m_N\}$}\label{subsub:glue21}

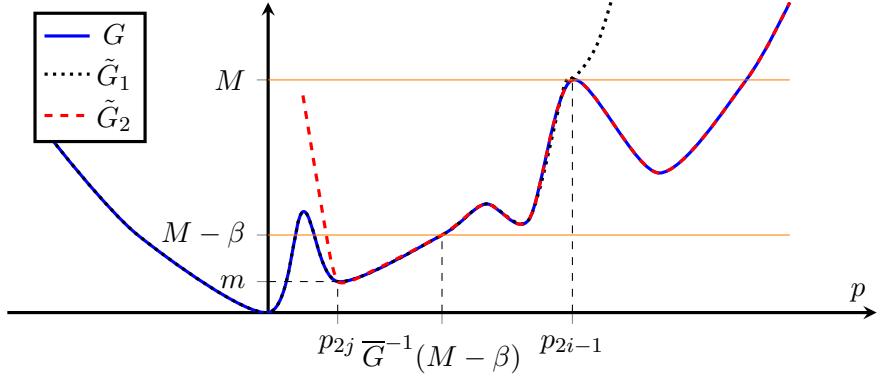
\begin{figure}[h!]
  \begin{tikzpicture}
  \begin{axis}[scale only axis=true,
        width=0.7\textwidth,
        height=0.25\textwidth,
        axis x line=middle, xlabel={$p$}, axis y line=middle,
    ylabel={$\ $}, tick align=outside, samples=50,
    xtick={0.8,2,3.5},
    xticklabels={$p_{2j}$,$\small\overline{G}^{-1}(M-\beta)$,$p_{2i-1}$}, ytick={0.4,1,3}, yticklabels={$m$,$M-\beta$,$M$}, xmin=-3, xmax=7, ymin=0, ymax=4, very thick, legend pos = north west]
    \addplot[smooth, no marks, blue, solid] coordinates {
      (-2.5,2.25)
      (-1.5,1)
      (0,0)
      (0.4,1.3)
      (0.8,0.4)
      (2,1)
      (2.5,1.4)
      (3,1.2)
      (3.5,3)
      (4.5,1.8)
      (5.5,3)
      (6,4)
    };
    \addlegendentry{$G$}
    \addplot+[smooth, no marks, black, dotted] coordinates {
      (-2.5,2.25)
      (-1.5,1)
      (0,0)
      (0.4,1.3)
      (0.8,0.4)
      (2,1)
      (2.5,1.4)
      (3,1.2)
      (3.4,2.9)
      (3.5,3)
    };
    \addlegendentry{$\tilde{G}_1$}
    
    \addplot+[smooth, no marks, red, dashed] coordinates {
      (0.8,0.4)
      (1,0.43)
      (2,1)
      (2.5,1.4)
      (3,1.2)
      (3.5,3)
      (4.5,1.8)
      (5.5,3)
      (6,4)
    };
    \addlegendentry{$\tilde{G}_2$}
    \addplot[smooth, no marks, red, dashed] coordinates {
      (0.4,2.8)
      (0.7,0.8)
      (0.8,0.4)
    };
    \addplot+[smooth, no marks, black, dotted] coordinates {
    (3.5,3)
    (3.8,3.5)
    (4.2,5)
    };  
    
    \addplot[smooth, no marks, orange, solid, thin] coordinates {
      (0,1)
      (6,1)
    };
    \addplot[smooth, no marks, orange, solid, thin] coordinates {
      (0,3)
      (6,3)
    };
    \addplot[smooth, no marks, black, dashed, thin] coordinates {
      (0,0.4)
      (0.8,0.4)
    };
    \addplot[smooth, no marks, black, dashed, thin] coordinates {
      (0.8,0.0)
      (0.8,0.4)
    };
    \addplot[smooth, no marks, black, dashed, thin] coordinates {
      (3.5,0.0)
      (3.5,3)
    };
    \addplot[smooth, no marks, black, dashed, thin] coordinates {
      (2,0.0)
      (2,1)
    };
\end{axis}
\end{tikzpicture}
\caption{Case $M-m\ge\beta>M-\min\{m_i,\ldots,m_N\}$.}
\label{fig2}
\end{figure}

In particular, $i>j$ in \eqref{eq:ij}.  Define
\[
\tilde G_1(p) = \begin{cases} G(p) &\text{if}\ p\leqslant p_{2i - 1},\\ I(p) &\text{if}\ p\geqslant p_{2i - 1},\end{cases}
\]
and
\[
\tilde G_2(p) = \begin{cases} I(p) &\text{if}\ p\leqslant p_{2j},\\G(p) &\text{if}\ p\geqslant p_{2j},\end{cases}
\]
where $I:\R\to[m,+\infty)$ is chosen such that it is strictly
decreasing and strictly increasing on $(-\infty,p_{2j}]$ and
$[p_{2i - 1},+\infty)$, respectively, $\tilde G_1, \tilde G_2\in\Ham$
and $G = \tilde G_1\wedge \tilde G_2$ (see Figure ~\ref{fig2}).

We will prove that, if $\tilde G_1,\tilde G_2\in\Hamone(a,V,\beta)$, then
$G\in\Hamone(a,V,\beta)$. Let
\[
\ol G^{-1}(M-\beta) = \max\{p\in\R:\, G(p) = M - \beta \}.
\]

Note that $\ol G^{-1}(M-\beta) \geqslant p_{2j}$ and $[\ol G^{-1}(M-\beta),p_{2i-1}] \in \mathcal{I}_M(G)$. We apply Lemma \ref{lem:kargo} and deduce that, for every $\omega\in\Omega$, $\Sol(M,[\ol G^{-1}(M-\beta),p_{2i-1}],\omega) \ne \emptyset$. Let
\begin{align*}
	&\ul f^M(x,\omega):=\inf\{f(x):\,f\in\Sol(M,[\ol G^{-1}(M-\beta),p_{2i-1}],\omega)\} \quad\text{and}\\
	&\ol f^M(x,\omega):=\sup\{f(x):\,f\in\Sol(M,[\ol G^{-1}(M-\beta),p_{2i-1}],\omega)\}
\end{align*}
as in \eqref{eq:lazim}. By Lemma \ref{lem:measurable}, these functions are jointly measurable and stationary. Moreover, $\ul f^M(\,\cdot\,,\omega), \ol f^M(\,\cdot\,,\omega)\in\Sol(M,[\ol G^{-1}(M-\beta),p_{2i-1}],\omega)$ for every $\omega\in\Omega$. By Lemma \ref{lem:cerca1}, for every $\delta > 0$, $y_0 > 0$ and $\P$-a.e.\ $\omega$, there is an interval $[L_1,L_2]$ such that
\[ \int_{L_1}^{L_2}\frac{dx}{a(x,\omega)} = y_0
\quad\text{and}\quad
\ol f^M(x,\omega) - \ul f^M(x,\omega) \leqslant \delta\ \ \forall x\in[L_1,L_2].
\]

Let
\[ \ul\theta(M) = \mathbb{E}[\ul f^M(0,\omega)]\in[p_{2j},p_{2i - 1}] \quad\text{and}\quad \ol\theta(M) = \mathbb{E}[\ol f^M(0,\omega)]\in[p_{2j},p_{2i - 1}]. \]
Fix any $\theta < \ul\theta(M)$. Note that $\tilde G_1\in\Hamzero\cap\Hamone(a,V,\beta)$ by our assumptions. Take $\lambda = \lambda(\theta)\geqslant \beta$ and $\ul f,\ol f$ as in Definition \ref{def:Hamone} for $\tilde G_1$. Since $G = \tilde G_1$ on $(-\infty,p_{2i-1}]$, it follows from the lemma below that $\ul f,\ol f$ are stationary solutions  to \eqref{eq:cellODE} with $G$.

\begin{lemma}\label{lem:bound3}
  With the notation in the paragraph above, we have the following
  dichotomy:
	\begin{itemize}
        \item [(i)]
          $\ul f(\,\cdot\,,\omega) \leqslant \ol f(\,\cdot\,,\omega)
          \leqslant 0 $ on $\mathbb{R}$ for $\P$-a.e.\ $\omega$;
		\item [(ii)] $\lambda < M$ and $0 \leqslant \ul f(\,\cdot\,,\omega) \leqslant \ol f(\,\cdot\,,\omega) < \ul f^M(\,\cdot\,,\omega) \leqslant p_{2i-1}$ on $\mathbb{R}$ for $\P$-a.e.\ $\omega$.
 	\end{itemize}
\end{lemma}

\begin{proof}
	By Lemma \ref{lem:clave} (with $\tilde G_1$), we are either in case (i) above or $0 \leqslant \ul f(\,\cdot\,,\omega) \leqslant \ol f(\,\cdot\,,\omega)$ on $\mathbb{R}$ for $\P$-a.e.\ $\omega$. Let us consider the latter case. We first show that $\lambda\ne M$. Suppose (for the purpose of reaching a contradiction) that 
$\lambda = M$. Again by Lemma \ref{lem:clave} (with $\tilde G_1$), for $\P$-a.e.\ $\omega$, $\ul f(\,\cdot\,,\omega)$ stays in an interval in $\mathcal{I}_M(\tilde G_1)$. There is a unique such interval that is a subset of $[0,+\infty)$, namely $[\ol G^{-1}(M-\beta),p_{2i-1}]$.
	Since $G = \tilde G_1$ on $(-\infty,p_{2i-1}]$ and $\ul f^M(\,\cdot\,,\omega)$ is minimal in $\Sol(M,[\ol G^{-1}(M-\beta),p_{2i-1}],\omega)$, it follows that $\ul f(\,\cdot\,,\omega) \geqslant \ul f^M(\,\cdot\,,\omega)$ on $\R$ for $\P$-a.e.\ $\omega$. This implies that $\theta \geqslant \mathbb{E}[\ul f(0,\omega)] \geqslant \mathbb{E}[\ul f^M(0,\omega)] = \ul\theta(M) > \theta$, which is a contradiction. 
	
Now that we know that $\lambda \neq M$, we can apply Corollary \ref{cor:ordqu} with $H(p,x,\omega) = G(p) + \beta V(x,\omega)$   
to infer that, for $\P$-a.e.\ $\omega$, the distance between $\ul f^M(\,\cdot\,,\omega)$ and $\ul f(\,\cdot\,,\omega)$ (resp.\ $\ol f(\,\cdot\,,\omega)$) is greater than $\delta|\lambda - M|$ for some $\delta>0$. Since $\mathbb{E}[\ul f^M(0,\omega)] = \ul\theta(M) > \theta \geqslant \mathbb{E}[\ul f(0,\omega)]$, we see that $\ul f(\,\cdot\,,\omega) + \delta|\lambda - M| \leqslant \ul f^M(\,\cdot\,,\omega)$ on $\mathbb{R}$ for $\P$-a.e.\ $\omega$. Recalling from Definition \ref{def:Hamone} that $\inf\{\ol f(x,\omega) - \ul f(x,\omega):\,x\in\R\} = 0$ for $\P$-a.e.\ $\omega$, we conclude that $\ol f(\,\cdot\,,\omega) + \delta|\lambda - M| \leqslant \ul f^M(\,\cdot\,,\omega)$ on $\mathbb{R}$ for $\P$-a.e.\ $\omega$.
	
	It remains to argue that $\lambda < M$. By Lemma \ref{lem:clave}, for $\P$-a.e.\ $\omega$, $\ul f(\,\cdot\,,\omega)$ and $\ol f(\,\cdot\,,\omega)$ stay in an interval in $\mathcal{I}_\lambda(G)$. However, no such interval is a subset of $[0,p_{2i-1}]$ if $\lambda > M$.
\end{proof}

Similarly, fix any $\theta > \ol\theta(M)$. Note that $\tilde G_2\in\Hamone(a,V,\beta)$ and $\tilde G_2(\,\cdot\,+ p_{2j}) - m \in\Hamzero$ by our assumptions. (Recall \eqref{eq:ij}). Take $\lambda = \lambda(\theta)\geqslant m + \beta$ and $\ul f,\ol f$ as in Definition \ref{def:Hamone} for $\tilde G_2$. Since $G = \tilde G_2$ on $[p_{2j},+\infty)$, it follows from the lemma below that $\ul f,\ol f$ are stationary solutions  to \eqref{eq:cellODE} with $G$.

\begin{lemma}\label{lem:bound4}
	With the notation in the paragraph above, $\lambda > M$ and $\ol f(\,\cdot\,,\omega) \geqslant \ul f(\,\cdot\,,\omega) > \ol f^M(\,\cdot\,,\omega) \geqslant \ol G^{-1}(M-\beta) \geqslant p_{2j}$ on $\mathbb{R}$ for $\P$-a.e.\ $\omega$.
\end{lemma}

\begin{proof}
	By Lemma \ref{lem:clave} (with $\tilde G_2$), we have the following dichotomy:
	\begin{itemize}
		\item [(i)] $\ul f(\,\cdot\,,\omega) \leqslant \ol f(\,\cdot\,,\omega) \leqslant p_{2j}$ on $\mathbb{R}$ for $\P$-a.e.\ $\omega$;
		\item [(ii)] $\ol f(\,\cdot\,,\omega) \geqslant \ul f(\,\cdot\,,\omega) \geqslant p_{2j}$ on $\mathbb{R}$ for $\P$-a.e.\ $\omega$.
	\end{itemize}
	Since $p_{2j} \leqslant \ol\theta(M) = \mathbb{E}[\ol f^M(0,\omega)] < \theta \leqslant \mathbb{E}[\ol f(0,\omega)]$, we rule out the first case.
	
	Let us first prove that $\lambda\ne M$. Suppose (for the purpose of reaching a contradiction) that $\lambda = M$. Again by Lemma \ref{lem:clave} (with $\tilde G_2$), for $\P$-a.e.\ $\omega$, $\ol f(\,\cdot\,,\omega)$ stays in an interval in $\mathcal{I}_M(\tilde G_2)$. There is a unique such interval that is a subset of $[p_{2j},+\infty)$, namely $[\ol G^{-1}(M-\beta),p_{2i-1}]$.
	Since $G = \tilde G_2$ on $[p_{2j},+\infty)$ and $\ol f^M(\,\cdot\,,\omega)$ is maximal in $\Sol(M,[\ol G^{-1}(M-\beta),p_{2i-1}],\omega)$, it follows that $\ol f(\,\cdot\,,\omega) \leqslant \ol f^M(\,\cdot\,,\omega)$ on $\R$ for $\P$-a.e.\ $\omega$. This implies that $\theta \leqslant \mathbb{E}[\ol f(0,\omega)] \leqslant \mathbb{E}[\ol f^M(0,\omega)] = \ol\theta(M) < \theta$, which is a contradiction.
		
Now that we know that $\lambda \neq M$, we can apply Corollary \ref{cor:ordqu} with $H(p,x,\omega) = G(p) + \beta V(x,\omega)$ to  infer that, for $\P$-a.e.\ $\omega$, the distance between  $\ol f^M(\,\cdot\,,\omega)$ and $\ol f(\,\cdot\,,\omega)$ (resp.\ $\ul f(\,\cdot\,,\omega)$) is greater than $\delta|\lambda - M|$ for some $\delta>0$. Since $\mathbb{E}[\ol f^M(0,\omega)] = \ol\theta(M) < \theta \leqslant \mathbb{E}[\ol f(0,\omega)]$, we see that $\ol f^M(\,\cdot\,,\omega) + \delta|\lambda - M| \leqslant \ol f(\,\cdot\,,\omega)$ on $\mathbb{R}$ for $\P$-a.e.\ $\omega$. Recalling from Definition \ref{def:Hamone} that $\inf\{\ol f(x,\omega) - \ul f(x,\omega):\,x\in\R\} = 0$ for $\P$-a.e.\ $\omega$, we conclude that $\ol f^M(\,\cdot\,,\omega) + \delta|\lambda - M| \leqslant \ul f(\,\cdot\,,\omega)$ on $\mathbb{R}$ for $\P$-a.e.\ $\omega$.
	
	It remains to argue that $\lambda > M$. By Lemma \ref{lem:clave}, for $\P$-a.e.\ $\omega$, $\ul f(\,\cdot\,,\omega)$ and $\ol f(\,\cdot\,,\omega)$ stay in an interval in $\mathcal{I}_\lambda(G)$. However, no such interval is a subset of $[\ol G^{-1}(M-\beta),+\infty)$ if $\lambda < M$.
\end{proof}

We recall the notation in Definition \ref{def:Hamone} and conclude that $G\in\Hamone(a,V,\beta)$ with the following choices:
\begin{itemize}
	\item if $\theta\in(-\infty,\ul\theta(M))$, then $\lambda \geqslant \beta$ and $\ul f,\ol f$ are as in Definition \ref{def:Hamone} for $\tilde G_1$;
	
	\item if $\theta\in[\ul\theta(M),\ol\theta(M)]$, then $\lambda = M$, $\ul f = \ul f^M$ and $\ol f = \ol f^M$;
	
	\item if $\theta\in(\ol\theta(M),+\infty)$, then $\lambda > M$ and $\ul f,\ol f$ are as in Definition \ref{def:Hamone} for $\tilde G_2$.
\end{itemize}

\subsubsection{$\beta \leqslant M - \min\{m_i,\ldots,m_N\}$}\label{subsub:glue22}

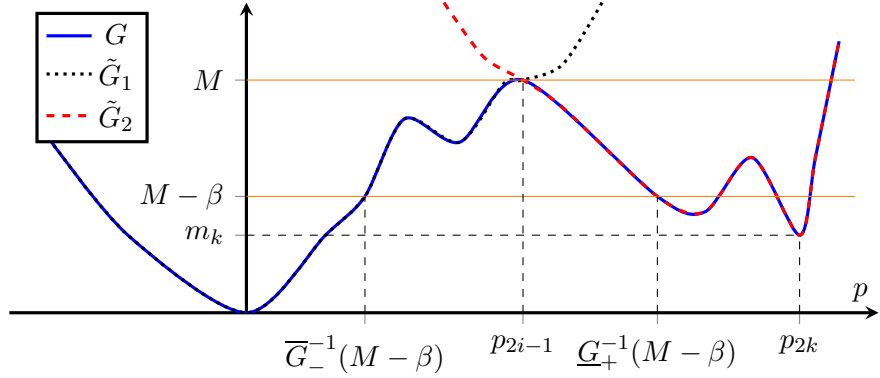
\begin{figure}[h!]
  \begin{tikzpicture}
  \begin{axis}[scale only axis=true,
        width=0.7\textwidth,
        height=0.25\textwidth,
        axis x line=middle, xlabel={$p$}, axis y line=middle,
    ylabel={$\ $}, tick align=outside, samples=100,
    xtick={1.5,3.5,5.2,7},
    xticklabels={$\overline{G}_-^{-1}(M-\beta)$,$p_{2i-1}$,$\underline{G}_+^{-1}(M-\beta)$,$p_{2k}$}, ytick={1,1.5,3}, yticklabels={$m_k$,$M-\beta$, $M$}, xmin=-3, xmax=8, ymin=0, ymax=4, very thick,legend pos = north west]
    \addplot[smooth, no marks, blue, solid] coordinates {
      (-2.5,2.25)
      (-1.5,1)
      (0,0)
      (1,1)
      (1.5,1.5)
      (2,2.5)
      (2.7,2.2)
      (3.5,3)
      (5.2,1.5)
      (5.8,1.3)
      (6.4,2) 
      (7,1)
      (7.2,2)
      (7.5,3.5)
    };
    \addlegendentry{$G$}
    \addplot+[smooth, no marks, black, dotted] coordinates {
      (-2.5,2.25)
      (-1.5,1)
      (0,0)
      (1,1)
      (1.5,1.5)
      (2,2.5)
      (2.7,2.2)
      (3.25,2.95)
      (3.5,3)
      (4,3.2)
      (4.5,4)
    };
    \addlegendentry{$\tilde{G}_1$}
    
    \addplot+[smooth, no marks, red, dashed] coordinates {
      (2.5,4)
      (3,3.3)
      (3.5,3)
      (4,2.65)
      (5.2,1.5)
      (5.8,1.3)
      (6.4,2) 
      (7,1)
      (7.2,2)
      (7.5,3.5)
    };
    \addlegendentry{$\tilde{G}_2$}
    \addplot[smooth, no marks, orange, solid, thin] coordinates {
      (0,1.5)
      (7.7,1.5)
    };
    \addplot[smooth, no marks, orange, solid, thin] coordinates {
      (0,3)
      (7.7,3)
    };
    \addplot[smooth, no marks, black, dashed, thin] coordinates {
      (0,1)
      (7,1)
    };
    \addplot[smooth, no marks, black, dashed, thin] coordinates {
      (1.5,0)
      (1.5,1.5)
    };
    \addplot[smooth, no marks, black, dashed, thin] coordinates {
      (3.5,0)
      (3.5,3)
    };
    \addplot[smooth, no marks, black, dashed, thin] coordinates {
      (5.2,0)
      (5.2,1.5)
    };
    \addplot[smooth, no marks, black, dashed, thin] coordinates {
      (7,0)
      (7,1)
    };
\end{axis}
\end{tikzpicture}
\caption{Case $\beta\le M-\min\{m_i,\ldots,m_N\}$.}
\label{fig3}
\end{figure}

By Definition \ref{def:Hamtwo}, in addition to \eqref{eq:ij},
\begin{equation}\label{eq:k}
	\text{there is a unique $k\in\{i,\ldots,N\}$ such that $m_k = \min\{m_i,\ldots,m_N\}$.}
\end{equation}
Define
\[
\tilde G_1(p) = \begin{cases} G(p) &\text{if}\ p\leqslant p_{2i - 1},\\ I(p) &\text{if}\ p\geqslant p_{2i - 1},\end{cases}
\]
and
\[
\tilde G_2(p) = \begin{cases} I(p) &\text{if}\ p\leqslant p_{2i-1},\\G(p) &\text{if}\ p\geqslant p_{2i-1},\end{cases}
\]
where $I:\R\to[M,+\infty)$ is chosen such that it is strictly
decreasing and strictly increasing on $(-\infty,p_{2i-1}]$ and
$[p_{2i - 1},+\infty)$, respectively, $\tilde G_1, \tilde G_2\in\Ham$
and $G = \tilde G_1\wedge \tilde G_2$ (see Figure~\ref{fig3}).

We will prove that, if $\tilde G_1,\tilde G_2\in\Hamone(a,V,\beta)$, then $G\in\Hamone(a,V,\beta)$. Let
\begin{align*}
	\ol G_-^{-1}(M-\beta) &= \max\{p\leqslant p_{2i-1}:\, G(p) = M-\beta \} \quad\text{and}\\
	\ul G_+^{-1}(M-\beta) &= \min\{p\geqslant p_{2i-1}:\, G(p) = M-\beta \}.
\end{align*}

Note that $[\ol G_-^{-1}(M-\beta),p_{2i-1}], [p_{2i-1},\ul G_+^{-1}(M-\beta)] \in \mathcal{I}_M(G)$. We apply Lemma \ref{lem:kargo} and deduce that, for every $\omega\in\Omega$, $\Sol(M,[\ol G_-^{-1}(M-\beta),p_{2i-1}],\omega) \ne \emptyset$ and $\Sol(M,[p_{2i-1},\ul G_+^{-1}(M-\beta)],\omega) \ne \emptyset$. Let
\begin{align*}
	&\ul f_-^M(x,\omega):=\inf\{f(x):\,f\in\Sol(M,[\ol G_-^{-1}(M-\beta),p_{2i-1}],\omega)\} \quad\text{and}\\
	&\ol f_+^M(x,\omega):=\sup\{f(x):\,f\in\Sol(M,[p_{2i-1},\ul G_+^{-1}(M-\beta)],\omega)\}
\end{align*}
as in \eqref{eq:lazim}. By Lemma \ref{lem:measurable}, these functions are jointly measurable and stationary. Moreover, $\ul f_-^M(\,\cdot\,,\omega)\in\Sol(M,[\ol G_-^{-1}(M-\beta),p_{2i-1}],\omega)$ and $\ol f_+^M(\,\cdot\,,\omega)\in\Sol(M,[p_{2i-1},\ul G_+^{-1}(M-\beta)],\omega)$ for every $\omega\in\Omega$. By Lemma \ref{lem:cerca2}, for every $\delta > 0$, $y_0 > 0$ and $\P$-a.e.\ $\omega$, there is an interval $[L_1,L_2]$ such that
\[ \int_{L_1}^{L_2}\frac{dx}{a(x,\omega)} = y_0
\quad\text{and}\quad
p_{2i-1} - \delta \leqslant \ul f_-^M(x,\omega) \leqslant p_{2i-1} \leqslant \ol f_+^M(x,\omega) \leqslant p_{2i-1} + \delta\ \ \forall x\in[L_1,L_2].
\]

Let
\begin{equation}
  \label{flat}
  \ul\theta_-(M) = \mathbb{E}[\ul f_-^M(0,\omega)] < p_{2i-1} \quad\text{and}\quad \ol\theta_+(M) = \mathbb{E}[\ol f_+^M(0,\omega)] > p_{2i-1}. 
\end{equation}
Fix any $\theta < \ul\theta_-(M)$. Note that $\tilde G_1\in\Hamzero\cap\Hamone(a,V,\beta)$ by our assumptions. Take $\lambda = \lambda(\theta)\geqslant \beta$ and $\ul f,\ol f$ as in Definition \ref{def:Hamone} for $\tilde G_1$. Since $G = \tilde G_1$ on $(-\infty,p_{2i-1}]$, it follows from the lemma below that $\ul f,\ol f$ are stationary solutions to \eqref{eq:cellODE} with $G$.

\begin{lemma}\label{lem:bound5}
	With the notation in the paragraph above, we have the following dichotomy:
	\begin{itemize}
        \item [(i)]
          $\ul f(\,\cdot\,,\omega) \leqslant \ol f(\,\cdot\,,\omega)
          \leqslant 0 $ on $\mathbb{R}$ for $\P$-a.e.\ $\omega$;
		\item [(ii)] $\lambda < M$ and $0 \leqslant \ul f(\,\cdot\,,\omega) \leqslant \ol f(\,\cdot\,,\omega) < \ul f_-^M(\,\cdot\,,\omega) \leqslant p_{2i-1}$ on $\mathbb{R}$ for $\P$-a.e.\ $\omega$.
	\end{itemize}
\end{lemma}

\begin{proof}
	Same as the proof of Lemma \ref{lem:bound3} after replacing $\ol G^{-1}(M-\beta)$, $\ul f^M$ and $\ul\theta(M)$ with $\ol G_-^{-1}(M-\beta)$, $\ul f_-^M$ and $\ul\theta_-(M)$, respectively.
\end{proof}

Similarly, fix any $\theta > \ol\theta_+(M)$. Note that $\tilde G_2\in\Hamone(a,V,\beta)$ and $\tilde G_2(\,\cdot\,+ p_{2k}) - m_k \in\Hamzero$ by our assumptions. (Recall \eqref{eq:k}). Take $\lambda = \lambda(\theta)\geqslant m_k + \beta$ and $\ul f,\ol f$ as in Definition \ref{def:Hamone} for $\tilde G_2$. Since $G = \tilde G_2$ on $[p_{2i-1},+\infty)$, it follows from the lemma below that $\ul f,\ol f$ are stationary solutions to \eqref{eq:cellODE} with $G$.

\begin{lemma}\label{lem:bound6}
	With the notation in the paragraph above, we have the following dichotomy:
	\begin{itemize}
		\item [(i)] $p_{2i-1} < p_{2k} \leqslant \ul f(\,\cdot\,,\omega) \leqslant \ol f(\,\cdot\,,\omega)$ on $\mathbb{R}$ for $\P$-a.e.\ $\omega$;
		\item [(ii)] $\lambda < M$ and $p_{2i-1} \leqslant \ol f_+^M(\,\cdot\,,\omega) < \ul f(\,\cdot\,,\omega) \leqslant \ol f(\,\cdot\,,\omega) \leqslant p_{2k}$ on $\mathbb{R}$ for $\P$-a.e.\ $\omega$.
	\end{itemize}
\end{lemma}

\begin{proof}
	Analogous to the proof of Lemma \ref{lem:bound5}.
\end{proof}

We recall the notation in Definition \ref{def:Hamone} and conclude that $G\in\Hamone(a,V,\beta)$ with the following choices:
\begin{itemize}
	\item if $\theta\in(-\infty,\ul\theta_-(M))$, then $\lambda \geqslant \beta$ and $\ul f,\ol f$ are as in Definition \ref{def:Hamone} for $\tilde G_1$;
	
	\item if $\theta\in[\ul\theta_-(M),\ol\theta_+(M)]$, then $\lambda = M$, $\ul f = \ul f_-^M$ and $\ol f = \ol f_+^M$; moreover, $\ul\theta_-(M)<\ol\theta_+(M)$ by \eqref{flat};
	
	\item if $\theta\in(\ol\theta_+(M),+\infty)$, then $\lambda \geqslant m_k + \beta$ and $\ul f,\ol f$ are as in Definition \ref{def:Hamone} for $\tilde G_2$.\\
\end{itemize}

\section{Proofs of Theorems \ref{thm:induction} and \ref{thm:main reduction}}\label{sec:mainproofs}

\subsection{Proof of Theorem \ref{thm:induction}}\label{sec:induction}

The proof is obtained by putting together our results in Sections \ref{sec:base} and \ref{sec:gluing}.
For every $N\geqslant 0$, let
\[ 
\HamzeroN = \{G\in\Hamzero:\, \text{$G$ has exactly $2N+1$ local
    extrema}\}. 
\] 
Note that $\Hamzero =\bigcup_{N=0}^\infty\HamzeroN$. Let us set $\Hamone:=\bigcap \Hamone(a,V,\beta)$, where the intersection is taken over all possible pairs of random functions $(a,V)$ satisfying (A), (V), (S) and $\beta>0$.  
We will prove by strong induction that
$\HamzeroN\subset\Hamone$ for every $N\geqslant0$:
\begin{itemize}
		
\item If $G\in\mathscr{G}_{0,0}$ (i.e., it is quasiconvex), then $G\in\Hamone$ by our result in Section \ref{sub:base1}.
	
\item Fix an $N\geqslant 1$. Suppose that $\mathscr{G}_{0,n}\subset\Hamone$ for all $n < N$. Take any $G\in\HamzeroN$. 
In view of the gluing result in Section \ref{sub:glue1}, it suffices to prove the assertion in the following two cases.\smallskip\\
{ {\em Case 1: $G$ is strictly decreasing on $(-\infty,0]$.}}\smallskip

Recall the notation in \eqref{eq:Mm}. We have two alternatives:\smallskip
\begin{itemize}
\item if $\beta > M - m$, then $G\in\Hamone$ by our result in Section \ref{sub:base2};\smallskip
\item if $\beta \leqslant M - m$, then recall \eqref{eq:ij}. We split the proof in the following two subcases.\smallskip
  \begin{itemize}
  \item [(i)] If $\beta > M - \min\{m_i,\ldots,m_N\}$, then we can
    write $G = \tilde G_1 \wedge \tilde G_2$ with
    $\tilde G_1, \tilde G_2\in\Ham$ given in Section
    \ref{subsub:glue21}. Note that
    $\tilde G_1(\,\cdot\,)\in\mathscr{G}_{0,n_1}$ and $\tilde G_2(\,\cdot\,+ p_{2j}) -
    m\in\mathscr{G}_{0,n_2}$ for some $n_1,n_2 < N$. So,
    $\tilde G_1(\,\cdot\,),\tilde G_2(\,\cdot\,+ p_{2j}) -
    m\in\Hamone$ by the induction hypothesis. Therefore, by an
    elementary change of variables (which is left to the reader),
    $\tilde G_2\in\Hamone$, too. We conclude that $G\in\Hamone$ by our
    result in Section \ref{subsub:glue21}.\smallskip
  \item [(ii)] If $\beta \leqslant M - \min\{m_i,\ldots,m_N\}$, then
    we can write $G = \tilde G_1 \wedge \tilde G_2$ with
    $\tilde G_1, \tilde G_2\in\Ham$ given in Section
    \ref{subsub:glue22}. Recall \eqref{eq:k}. Note that
    $\tilde G_1(\,\cdot\,)\in\mathscr{G}_{0,n_1}$ and $\tilde G_2(\,\cdot\,+ p_{2k}) -
    m_k\in\mathscr{G}_{0,n_2}$ for some $n_1,n_2 < N$. Therefore,
    $\tilde G_1(\,\cdot\,),\tilde G_2(\,\cdot\,+ p_{2k}) -
    m_k\in\Hamone$ by the induction hypothesis, and
    $\tilde G_2\in\Hamone$, too, as above. We conclude that
    $G\in\Hamone$ by our result in Section
    \ref{subsub:glue22}.\smallskip
  \end{itemize}
	\end{itemize}
{\em Case 2: $G$ is strictly increasing on $[0,+\infty)$.}\smallskip

Let us set $\check G(p):=G(-p)$ for all $p\in\R$. By Case 1, we have that $\check G\in\Hamone$. We conclude that $G\in\Hamone$ after 
noticing that  $G\in \Hamone(a,V,\beta)$ if and only if $\check G\in\Hamone(\check a,\check V,\beta)$, where $\check a(x,\omega):=a(-x,\omega)$ and $\check V(x,\omega):=V(-x,\omega)$. To see this, it suffices to notice that $f$ is a $\CC^1_b$ solution of equation \eqref{eq:cellODE} if and only if $\check f(x):=-f(-x),\ x\in\R$, is a $\CC^1_b$ solution of  \eqref{eq:cellODE} with $\check a, \check V$ in place of $a, V$.\qed	
\end{itemize}

\subsection{Proof of Theorem \ref{thm:main reduction}}\label{sec:stable} 
Pick $G\in\Ham$. By (G1)--(G2),  there exists a $p_{\rm{min}}\in\R$ such that 
$G(p_{\rm{min}})=m_0:=\inf\{G(p):\,p\in\R\}$. Let us set $\widetilde G(\,\cdot\,):=G(\,\cdot\, + p_{\rm{min}}) - m_0$. An easy computation shows that homogenization (of equation \eqref{eq:generalHJ}) holds for $G$ if and only if it holds for $\widetilde G$ (and, in this instance, $\HV(\widetilde G)(\theta)=\HV(G)(\theta+ p_{\rm{min}})-m_0$ for all $\theta\in\R$). Hence, without any loss of generality, it suffices to prove the assertion for a $G\in\Ham$ additionally satisfying $G(p)\geqslant G(0)=0$ for all $p\in\R$. Let $\alpha_0,\alpha_1>0$ and $\gamma>1$ be constants such that $G$ belongs to the family $\Ham(\alpha_0,\alpha_1,\gamma)$ of functions that satisfy  
(G1)--(G2) with the triple $(\alpha_0,\alpha_1,\gamma)$. 
For every fixed $n\in\N$, choose $G_n\in\Hamzero\cap\Ham(\alpha_0,\alpha_1,\gamma)$ such that 
$\|G-G_n\|_{L^\infty([-n,n])}<1/n$ and $G_n$ is strictly decreasing (resp., increasing) on $(-\infty,-n]$ (resp., $[n,+\infty)$). By assumption, equation \eqref{eq:generalHJ} homogenizes for each $G_n$. By invoking Theorem \ref{appB teo stability}, we conclude that  \eqref{eq:generalHJ} homogenizes for $G$ as well.\qed

\appendix

\section{Stationary sub and supersolutions}\label{app:general}

Some of our arguments in Sections \ref{sec:spehom}--\ref{sec:gluing} 
rely on several results that hold true for rather general continuous and coercive Hamiltonians, which we present here for future reference. 

Throughout this appendix, we will work with a general probability space $(\Omega,\F, \P)$, where $\P$ and $\F$ denote the probability measure on $\Omega$  and
the $\sigma$--algebra of $\P$--measurable subsets of $\Omega$, respectively. 
We will assume that  $\P$ is invariant under the action of a one-parameter group $(\tau_x)_{x\in\R}$ of transformations $\tau_x:\Omega\to\Omega$ and that the action of $(\tau_x)_{x\in\R}$ is ergodic. (See the beginning of Section \ref{sec:result} for details.) 
However, no topological or completeness assumptions are made on the probability space.

Let $a:\R\times\Omega\to [0,1]$ be a stationary process, and $H:\R\times\R\times\Omega\to[0,+\infty)$ a stationary process with respect to the shifts in the second variable. We will assume that $a$ is continuous in the first variable and $H$ is continuous in the first two variables, for any fixed $\omega\in\Omega$.

For $\lambda \geqslant 0$, we will establish suitable comparison results for sub and supersolutions of the equation
\begin{equation*} 
	a(x,\omega)f'(x,\omega) + H(f(x,\omega),x,\omega) = \lambda,\quad x\in\R,
\end{equation*}
and eventually find stationary $\CC^1$ solutions, at least under proper assumptions on $H$. 

We introduce the following conditions on $H$ and specify in each statement which of them are needed for that result:
\begin{itemize}
	\item[(H1)] \quad there exist two coercive functions
	$\alpha_H,\beta_H:[0,+\infty)\to\R$ such that
	\[
	\alpha_H\left(|p|\right)\leqslant H(p,x,\omega)\leqslant \beta_H\left(|p|\right)\quad\hbox{for every $(p,x,\omega)\in\R\times\R\times\Omega$;}
	\]
	\item[(H2)] \quad for every fixed $R>0$, there exists a constant $C_R>0$ such that 
	\[
	|H(p,x,\omega)-H(q,x,\omega)|\leqslant C_R |p-q|\quad\hbox{for every $p,q\in [-R,R]$ and $(x,\omega)\in\R\times\Omega$}.
	\] 
\end{itemize}

For our first statement, we will need the following notation: for each $\lambda\geqslant 0$,
\begin{align}\label{eq:pla} 
	p_\lambda^-:=\inf_{x\in\R}\inf\{p\in\R:\ H(p,x,\omega)\leqslant \lambda\} \quad\text{and}\quad p^+_\lambda:=\sup_{x\in\R}\sup\{p\in\R:\ H(p,x,\omega)\leqslant \lambda\}.
\end{align}
By the ergodicity assumption and (H1), the quantities $p^\pm_\lambda$ are a.s.\ constants. The functions $\lambda\mapsto p^-_\lambda$ and $\lambda\mapsto p^+_\lambda$ are, respectively, non-increasing and non-decreasing (and, in general, not continuous). Furthermore, $p^-_\lambda\to-\infty$ and $p^+_\lambda\to +\infty$ as $\lambda\to+\infty$.

\begin{lemma}\label{lem:ubounds}
	Assume that $H$ satisfies {\em (H1)}. Take any $\lambda>0$. Let $f(x,\omega)$ be a stationary process such that, for all $\omega\in\Omega$, $f(\,\cdot\,,\omega)\in\CC^1(\R)$, and
	\begin{equation}\label{eq:lesslam}
		a(x,\omega)f'(x,\omega)+H(f(x,\omega),x,\omega)< \lambda\quad\forall x\in\R. 
	\end{equation}
	Then, on a set $\Omega_f$ of probability 1, $f(x,\omega)\in (p^-_\lambda,p^+_\lambda)$ for all $x\in\R$, where $p^\pm_\lambda$ are defined in \eqref{eq:pla}.
\end{lemma}

\begin{proof}
	Let $E_\lambda^-(\omega) = \{x\in\R: f(x,\omega) \leqslant p^-_\lambda\}$. By the ergodicity assumption, this set is empty with probability $0$ or $1$. Suppose that $\P(E_\lambda^-(\omega)\ne\emptyset)=1$. Then, $\inf E_\lambda^-(\omega)=-\infty$ for $\P$-a.e.\ $\omega$. Fix such an $\omega$ and take a sequence $x_n\in E_\lambda^-(\omega)$, $n\in\N$, such that $x_n\to-\infty$ as $n\to +\infty$.
	
	For each $n\in\N$, let $\ol{x}_n:=\sup\{x > x_n:\ [x_n,x)\subset E_\lambda^-(\omega)\}$.  We claim that $\ol{x}_n = +\infty$. First of all, the set over which we are taking a supremum is nonempty. Indeed, by \eqref{eq:lesslam} and the fact that $f(x_n,\omega) \leqslant p^-_\lambda$, we have
	\[a(x_n,\omega)f'(x_n,\omega)<\lambda - H(f(x_n,\omega),x_n,\omega) \leqslant 0.\]
	This implies that $f'(x_n,\omega)<0$ and, therefore, there is a $\delta>0$ such that $f(x,\omega) \leqslant f(x_n,\omega) \leqslant p^-_\lambda$ for all $x\in[x_n,x_n+\delta)$. If $\ol{x}_n < +\infty$, then $f(\ol{x}_n,\omega) = p^-_\lambda$ and, by the same argument, we get a contradiction with the definition of $\ol{x}_n$. We conclude that $[x_n, +\infty)\subset E_\lambda^-(\omega)$ for all $n\in\N$ and, thus, $E_\lambda^-(\omega)=\R$. But then, with probability $1$, $H(f(x,\omega),x,\omega)\geqslant \lambda$ for all $x\in\R$ and
	\[a(x,\omega)f'(x,\omega) < \lambda - H(f(x,\omega),x,\omega)\leqslant 0\quad\forall x\in\R.\]
	We deduce that $f'(x,\omega)<0$ for all $x\in\R$, and $x\mapsto f(x,\omega)$ is strictly decreasing on $\R$. In light of Lemma \ref{lem:stat} below, this contradicts the assumption that $f(x,\omega)$ is stationary. We conclude that $\P(E_\lambda^-(\omega)=\emptyset)=1$.
	
	If we let $E_\lambda^+(\omega) = \{x\in\R:\ f(x,\omega)\geqslant p^+_\lambda\}$, then $\P(E_\lambda^+(\omega)=\emptyset)=1$ by a symmetric argument.\footnote{Equivalently, we can apply the previous argument to the process $\check f(x,\omega):=-f(-x,\omega)$.} We conclude that, for every $\omega$ from a set $\Omega_f$ of full measure, $f(x,\omega)\in (p^-_\lambda,p^+_\lambda)$ for all $x\in\R$.
\end{proof}

\begin{lemma}\label{lem:stat}
	Let $f(x,\omega)$ be a stationary process such that $f(\,\cdot\,,\omega) \in\CC(\R)$ for every $\omega\in\Omega$. Then, we have the following dichotomy:
	\begin{itemize}
		\item [(i)] $\P(f(x,\omega) = c\ \forall x\in\R) = 1$ for some constant $c\in\R$;
	    \item [(ii)] for $\P$-a.e.\ $\omega$, $f(\,\cdot\,,\omega)$ has infinitely many local maxima and minima.
	\end{itemize}
\end{lemma}

\begin{proof}
	If $\P(f(0,\omega) = c) = 1$ for some $c\in\R$, then we are in case (i) by stationarity and continuity. Otherwise, there exist constants $c_1 < c_2$ such that $\P(f(0,\omega) \leqslant c_1) > 0$ and $\P(f(0,\omega) \geqslant c_2) > 0$. By stationarity and ergodicity, for $\P$-a.e.\ $\omega$, there exist $x_n$, $y_n$, $n\in\Z$, such that $f(x_n) \leqslant c_1$, $f(y_n) \geqslant c_2$, $x_n < y_n < x_{n+1}$, and $x_n\to\pm\infty$ as $n\to\pm\infty$. Therefore, $f(\,\cdot\,,\omega)$ has a local maximum (resp.\ local minimum) in each interval $[x_n,x_{n+1}]$ (resp.\ $[y_n,y_{n+1}]$), and we are in case (ii).
\end{proof}

As a consequence of Lemma \ref{lem:ubounds}, we infer
\begin{cor}\label{cor:lelam}
	Assume that $H$ satisfies {\em (H1)}. Take any $\lambda\geqslant 0$. Let $f(x,\omega)$ be a stationary process such that, for all $\omega\in\Omega$, $f(\,\cdot\,,\omega)\in\CC^1(\R)$, and
	\[ a(x,\omega)f'(x,\omega)+H(f(x,\omega),x,\omega)\leqslant \lambda\quad\forall x\in\R. \]
	Then, on a set $\Omega_f$ of probability 1,
	$\displaystyle f(x,\omega)\in[\sup_{\mu>\lambda}p^-_\mu,\inf_{\mu>\lambda}p^+_\mu]$ for all $x\in\R$.
\end{cor}

The next result generalizes the fact that, under suitable conditions, two distinct solutions of an ODE do not touch each other.

\begin{lemma}\label{lem:order}
	Assume that $H$ satisfies {\em (H2)}, and $a(x,\omega) > 0$ for all $(x,\omega)\in\R\times\Omega$. Let $f_1(x,\omega)$ and $f_2(x,\omega)$ be stationary processes such that, for all $\omega\in\Omega$, $f_1(\,\cdot\,,\omega),f_2(\,\cdot\,,\omega)\in\CC^1(\R)\cap \CC_b(\R)$, and
	\begin{equation}\label{eq:ord1}
		a(x,\omega)f_1'(x,\omega)+H(f_1(x,\omega),x,\omega) \leqslant a(x,\omega)f_2'(x,\omega)+H(f_2(x,\omega),x,\omega)\quad \forall x\in\R.
	\end{equation}
	Then, one of the following events has probability $1$:
	\begin{align*}
		\Omega_0\, &= \{\omega\in\Omega:\, (f_1-f_2)(x,\omega) = 0\ \text{for all}\ x\in\R \};\\
		\Omega_- &= \{\omega\in\Omega:\, (f_1-f_2)(x,\omega) < 0\ \text{for all}\ x\in\R \};\\
		\Omega_+ &= \{\omega\in\Omega:\, (f_1-f_2)(x,\omega) > 0\ \text{for all}\ x\in\R \}.
	\end{align*}
\end{lemma}

\begin{proof}
	The events $\Omega_0, \Omega_-,\Omega_+$ as well as the events
	\begin{align*}
		\wh\Omega_- &= \{\omega\in\Omega:\, (f_1-f_2)(x_0,\omega) < 0\ \text{for some}\ x_0\in\R \}\quad\text{and}\\
		\wh\Omega_+ &= \{\omega\in\Omega:\, (f_1-f_2)(x_0,\omega) > 0\ \text{for some}\ x_0\in\R \}
	\end{align*}
	are translation invariant, so each of them has probability $0$ or $1$ by ergodicity. If $\P(\Omega_0) = 0$, then $\P(\wh\Omega_-) = 1$ or $\P(\wh\Omega_+) = 1$ since 
	$\Omega = \Omega_0\cup \wh\Omega_-\cup\wh\Omega_+$. Assume that $\P(\wh\Omega_-) = 1$.
	\footnote{Otherwise we can work with
	\[ \check  H(p,x) := H(-p,-x),\quad \check  a(x,\omega) := a(-x,\omega),\quad \check  f_1(x,\omega) := -f_1(-x,\omega)\quad\text{and}\quad \check  f_2(x,\omega) := -f_2(-x,\omega). \]}
	
	For every $\omega\in\wh\Omega_-$, let $x_0 = x_0(\omega)$ be as in the definition of $\wh\Omega_-$, and let
	\[ \ol{x}_0 = \ol{x}_0(\omega) = \sup\{x>x_0: (f_1-f_2)(y,\omega) < 0\ \text{for all}\ y\in[x_0,{x})\}. \]
	We claim that $\ol{x}_0 = +\infty$. Suppose (for the sake of reaching a contradiction) that $\ol{x}_0 < +\infty$. Then, $g(x,\omega) := (f_2-f_1)(x,\omega) > 0$ on $[x_0,\ol{x}_0)$, and $g(\ol{x}_0,\omega) = 0$. By \eqref{eq:ord1} and (H2),
	\begin{align*}
		a(x,\omega)g'(x,\omega) &= a(x,\omega)(f_2'(x,\omega) - f_1'(x,\omega))\\
		                        &\geqslant H(f_1(x,\omega),x,\omega) - H(f_2(x,\omega),x,\omega) \geqslant -C_Rg(x,\omega) \quad \forall x\in[x_0,\ol{x}_0],
	\end{align*}
	where $R = \max\{\|f_1\|_\infty,\|f_2\|_\infty\}$. Therefore,
	\[ \log g(\ol x_0) = \log g(x_0) + \int_{x_0}^{\ol x_0} \frac{g'(x,\omega)}{g(x,\omega)} dx \geqslant \log g(x_0) - C_R \int_{x_0}^{\ol x_0}\frac{dx}{a(x,\omega)} > -\infty, \]
	which is a contradiction. We deduce that $\ol{x}_0 = +\infty$ is true, i.e., $(f_1-f_2)(y,\omega) < 0$ for all $y \geqslant x_0$.
	
	For every $\omega\in\Omega$, let
	\[ \ul{x}_0 = \ul{x}_0(\omega) = \inf\{x\in\R:\,(f_1-f_2)(y,\omega) < 0\ \text{for all}\ y>x\}. \]
	Observe that, by stationarity, $\P(\ul{x}_0(\omega)\leqslant x)=\P(\ul{x}_0(\tau_z\omega)\leqslant x)=\P(\ul{x}_0(\omega)\leqslant x+z)$ for all $x,z\in\R$. This cannot hold for any $\R$-valued random variable, so we deduce that $\P(\ul{x}_0(\omega) \in\{-\infty,+\infty\})=1$. Since we have shown above that $\ul{x}_0(\omega) \leqslant x_0(\omega) < +\infty$ for every $\omega\in\wh\Omega_-$, we conclude that
	\[ \P(\Omega_-) = \P(\omega\in\Omega:\,\ul{x}_0(\omega) = -\infty)=1.\qedhere \]
\end{proof}

By strengthening the inequality in \eqref{eq:ord1}, we get a quantitative estimate on the difference between $f_1$ and $f_2$.

\begin{lemma}\label{lem:ordqu}
	Assume that $H$ satisfies {\em (H2)}. Take any $\epsilon>0$. Let $f_1(x,\omega)$ and $f_2(x,\omega)$ be stationary processes such that, for all $\omega\in\Omega$, $f_1(\,\cdot\,,\omega), f_2(\,\cdot\,,\omega)\in\CC^1(\R)\cap \CC_b(\R)$, and
	\begin{equation}
		\label{eq:ord3}
		a(x,\omega)f_1'(x,\omega) + H(f_1(x,\omega),x,\omega)+\epsilon<a(x,\omega)f_2'(x,\omega) + H(f_2(x,\omega),x,\omega)\quad\forall x\in\R.
	\end{equation}
	Then, 
	\[ \P((f_1-f_2)(x,\omega) > \epsilon/{C_R}\ \forall x\in\R)=1\ \ \text{or}\ \ \P((f_2-f_1)(x,\omega) > \epsilon/{C_R}\ \forall x\in\R) = 1, \]
	where $R=\max\{\|f_1\|_\infty,\|f_2\|_\infty\}$.
\end{lemma}

\begin{proof}
	As we argued in the proof of Lemma \ref{lem:order}, by ergodicity it is sufficient to show that, with probability $1$,
	$|f_1(x,\omega)-f_2(x,\omega)| > \epsilon/{C_R}$ for all $x\in\R$. 
	
	{\em Case 1.} Suppose that, 
	with probability $1$, $f_2(x,\omega)\equiv f_1(x,\omega)+c$ for some constant $c\in\R$. By (H2) and \eqref{eq:ord3}, we get
	\begin{equation*}
		C_R|c|\geqslant H(f_1(x,\omega)+c,x,\omega) - H(f_1(x,\omega),x,\omega)>\epsilon.
	\end{equation*}
	We conclude that, for $\P$-a.e.\ $\omega$, $|f_1(x,\omega)-f_2(x,\omega)|\equiv |c|>\epsilon/C_R$.
	
	{\em Case 2.} If we are not in Case 1, then, with probability $1$, 
	$f_1(\,\cdot\,,\omega) - f_2(\,\cdot\,,\omega)$ has infinitely many local maxima and minima by Lemma \ref{lem:stat}. By (H2) and \eqref{eq:ord3}, at each local extremum $x_m$,
	\[ C_R|f_1(x_m,\omega)-f_2(x_m,\omega)| \geqslant H(f_2(x_m,\omega),x_m,\omega) - H(f_1(x_m,\omega),x_m,\omega)>\epsilon. \]
	We conclude that, for $\P$-a.e.\ $\omega$, $|f_1(x,\omega)-f_2(x,\omega)| > \epsilon/C_R$ for all $x\in\R$.
\end{proof}

As a consequence of the above results, the following statement holds.

\begin{corollary}\label{cor:ordqu}
	Assume that $H$ satisfies {\em (H1)} and {\em (H2)}. Take any $\lambda_1,\lambda_2\in\R$ such that $0\leqslant\lambda_1 < \lambda_2$. Let $f_1(x,\omega)$ and $f_2(x,\omega)$ be stationary processes such that, for all $i\in\{1,2\}$ and $\omega\in\Omega$, $f_i(\,\cdot\,,\omega)\in\CC^1(\R)$, and
	\begin{equation*}
		a(x,\omega)f_i'(x,\omega) + H(f_i(x,\omega),x,\omega) = \lambda_i\quad\forall x\in\R.
	\end{equation*}
	Then, there is a constant $\delta>0$, which depends only on $\lambda_2$ and $H$, such that
	\[\P((f_1-f_2)(x,\omega)>\delta(\lambda_2-\lambda_1)\ \forall x\in\R)=1\ \ \text{or}\ \ \P((f_2-f_1)(x,\omega)>\delta(\lambda_2-\lambda_1)\ \forall x\in\R)= 1.\]
\end{corollary}

\begin{proof}
	By Corollary \ref{cor:lelam}, on a set $\Omega_{f_1}\cap\Omega_{f_2}$ of probability $1$,
	$f_i(x,\omega)\in\displaystyle{[\sup_{\mu>\lambda_2}p^-_\mu,\inf_{\mu>\lambda_2}p^+_\mu]}$ for all $i\in\{1,2\}$ and $x\in\R$. Applying Lemma \ref{lem:ordqu} (with $\epsilon  = \lambda_2 - \lambda_1$ and after restricting the processes $f_1,f_2$ to $\Omega_{f_1}\cap\Omega_{f_2}$), we deduce that the desired result holds with $\delta  = 1/{C_R}$, where $R = \displaystyle{\max\{ |\sup_{\mu>\lambda_2}p^-_\mu|,|\inf_{\mu>\lambda_2}p^+_\mu| \}}$.
\end{proof}

The remaining statements are deterministic. One should think of $a(\,\cdot\,)$ and $H(\,\cdot\,,\,\cdot\,)$ as $a(\,\cdot\,,\omega)$ and $H(\,\cdot\,,\,\cdot\,,\omega)$ with $\omega$ fixed.

\begin{lemma}\label{lem:charac}
	Assume that $H(p,x)$ satisfies {\em (H2)}, and $a(x)>0$ for all $x\in\R$. Let $f_1, f_2\in\CC^1(\R)\cap \CC_b(\R)$ be such that $\inf\{f_2(x) - f_1(x):\,x\in\R\} = 0$ and
	\[ a(x)f_1'(x) + H(f_1(x),x) = a(x)f_2'(x) + H(f_2(x),x),\quad x\in\R. \]
	Then, for every $\delta > 0$ and $y_0 > 0$, there is an interval $[L_1,L_2]$ such that
	\[ \int_{L_1}^{L_2}\frac{dx}{a(x)} = y_0
	\quad\text{and}\quad
	f_2(x) - f_1(x) \leqslant \delta,\ \ \forall x\in[L_1,L_2].
	\]
\end{lemma}

\begin{proof}
  Assume that $g(x) := f_2(x) - f_1(x)$ is not
  identically zero. (Otherwise, the desired conclusion trivially
  holds.) Note that 
	\begin{equation}\label{eq:ssb}
		a(x)|g'(x)| = a(x)|f'_2(x) - f'_1(x)| = |H(f_1(x),x) - H(f_2(x),x)| \leqslant C_R g(x),\quad\forall x\in\R,
	\end{equation}
	by (H2), where $R = \max\{\|f_1\|_\infty,\|f_2\|_\infty\}$. For every $\delta > 0$ and $y_0 > 0$, there exist $L_1,L_2\in\R$ with $L_1 < L_2$ such that $0 < g(L_1) \leqslant \delta e^{-C_Ry_0}$ and $\displaystyle{ \int_{L_1}^{L_2}\frac{dx}{a(x)} = y_0 }$\footnote{Since $a(x)\in(0,1]$,  $L-\ell\le \displaystyle \int_{\ell}^L\frac{dx}{a(x)}<+\infty$ for all $\ell,\,L\in\mathbb{R}$.}. For every $x\in [L_1,L_2]$,
	\[ \left|\log\left(\frac{g(x)}{g(L_1)}\right)\right| = \left|\int_{L_1}^{x} \frac{g'(y)}{g(y)}dy\right| \leqslant C_R \int_{L_1}^{x}\frac{dy}{a(y)} \leqslant C_Ry_0 \]
	by \eqref{eq:ssb} and, therefore, $0 < g(x) \leqslant \delta$.
\end{proof}

The next lemma states that we can always insert a global $\CC^1$ solution between two 
strict sub and supersolutions which do not intersect. The result is standard, we give a proof for the reader's convenience. 

\begin{lemma}\label{lem:inbetw}
	Assume that $H(p,x)$ satisfies {\em (H1)} and {\em (H2)}, $a(x)>0$ for all $x\in\R$, there exist bounded functions $m,M\in\CC^1(\R)$ such that $m(x)<M(x)$ for all $x\in\R$, and either one of the following hold:
	\begin{align}
		&a(x)m'(x)+H(m(x),x) < 0 < a(x)M'(x)+H(M(x),x)\quad\forall x\in\R,\ \text{or}\label{eq:catla}\\
		&a(x)m'(x)+H(m(x),x) > 0 > a(x)M'(x)+H(M(x),x)\quad\forall x\in\R.\label{eq:patla}
	\end{align}
	Then, there exists a function $ f\in\CC^1(\R)$ that solves the equation
	\begin{align}\label{eq:zero}
		a(x)f'(x)+H(f(x),x)=0\quad &\forall x\in \R, 
		\shortintertext{and satisfies}
		m(x) < f(x) < M(x)\quad&\forall x\in\R.\nonumber
	\end{align}
\end{lemma}

\begin{proof}
	Suppose that \eqref{eq:catla} holds.\footnote{If \eqref{eq:patla} holds, then we can work with the following for which \eqref{eq:catla} holds:
	\[ \check  H(p,x) := H(-p,-x),\quad \check  a(x,\omega) := a(-x,\omega),\quad \check  m(x) := -M(-x)\quad\text{and}\quad \check  M(x) := -m(-x). \]}
	For every $n\in \N$, consider a solution $f_n(x)$ of \eqref{eq:zero} on $[-n,+\infty)$ with the initial condition $f_n(-n) = m(-n)$. First of all, we note that such a solution exists, is unique, and satisfies $m(x)< f_n(x)<M(x)$ for all $x\in(-n,+\infty)$.  Indeed, by our conditions on the coefficients, we have local existence and uniqueness. The absence of a blow up and the claimed bounds can be seen as follows.
	
	Suppose (for the sake of reaching a contradiction) that the set $E:=\{x> -n:\,f_n(x)\leqslant m(x)\}$ is nonempty and set   
	$\ul{x}:=\inf E$. If  $\ul{x}>-n$, then $f_n(\ul{x})=m(\ul{x})$, $f_n'(\ul{x})>m'(\ul{x})$ by \eqref{eq:catla} and \eqref{eq:zero}, and, hence, $\exists\delta>0$ such that
	$f_n(x)<m(x)$ for all $x\in(\ul{x}-\delta,\ul{x})$, which is a contradiction with the definition of $\ul{x}$. If $\ul{x}=-n$, then by the same argument, $\exists\delta>0$ such that $f_n(x)>m(x)$ for all $x\in(\ul{x},\ul{x}+\delta)$. This contradicts again the definition of $\ul{x}$. We conclude that $f_n(x)>m(x)$ for all $x> -n$. The upper bound $f_n(x)<M(x)$ for all $x>-n$ is obtained in a similar way. 
	
	Fix $L\in\N$. Restricted to the interval $[-L,L]$, the sequence $(f_n)_{n\geqslant L}$ is monotone increasing in $n$, uniformly bounded and, since each $f_n$ is a solution of \eqref{eq:zero} with $a>0$, it is also equi-continuous. Therefore, the sequence $(f_n)_{n\geqslant L}$ converges uniformly on $[-L,L]$ as $n\to+\infty$ to a continuous function $f$. Since $f$ is a uniform limit of solutions $f_n$ of \eqref{eq:zero} on $[-L,L]$ and $a>0$, it is a bounded $\CC^1$ solution of \eqref{eq:zero} on $(-L,L)$. Letting $L\to+\infty$, we extend $f$ as a bounded $\CC^1$ solution on $\R$.
	
	By construction, $m(x) < f(x) \leqslant M(x)$ for all $x\in\R$. In fact, the second inequality is strict, too. Indeed, if $f(\ol x) = M(\ol x)$ for some $\ol x\in\R$, then 
	$f'(\ol{x}) < M'(\ol{x})$ by \eqref{eq:patla} and \eqref{eq:zero}, and, hence, $\exists\delta>0$ such that $f(x) > M(x)$ for all $x\in(\ol{x}-\delta,\ol{x})$, which is a contradiction.
\end{proof}
We end this section with the following useful stability result. 
\begin{lemma}\label{appendix A lemma lattice}
Assume $H(p,x)$ satisfies (H2) and $a(x)>0$ for all $x\in\R$. Let $\Sol_\lambda$ be a nonempty family of solutions to
\begin{equation}\label{appendix A eq ODE}
	a(x,\omega)u'(x) + H(u(x),x) = \lambda,\quad x\in\R.
\end{equation}
If $\Sol_\lambda$ is a compact subset of $\D{C}(\R)$, then the functions 
\[
\ul u (x):=\inf_{u\in\Sol_\lambda} u(x),
\quad
\ol u(x):=\sup_{u\in\Sol_\lambda} u(x),
\quad
x\in\R,
\] 
are in $\Sol_\lambda$. 
\end{lemma}

\begin{proof}
We only prove the assertion for $\ul u$, as the other case is analogous. Let $(u_n)_{n\in\N}\subset \Sol_\lambda$ such that $u_n(0)\searrow \ul u(0)$.  Since the graphs of the functions $u_n$ cannot cross, we have 
\[
u_n(x)>u_{n+1}(x)\quad\hbox{for all $n\in\N$ and $x\in\R$.}
\] 
Let us set $v(x):=\lim_{n} u_n(x)$ for all $x\in\R$. Since $\Sol_\lambda$ is compact in $\CC(\R)$, this limit is actually locally uniform in $\R$, and $v\in\CC(\R)$. By \eqref{appendix A eq ODE}, the derivatives $(u_n')_{n\in\N}$ also converge in $\CC(\R)$, hence the functions $u_n$  actually converge to $v$ in the local $\CC^1$ topology, yielding  $v\in\Sol_\lambda$. Furthermore, $\ul u(0)=v(0)$, hence $\ul u=v\in\Sol_\lambda$, as it was asserted. 
\end{proof}

\section{PDE material}\label{app:PDE}

In this appendix, we collect some PDE results that we need in the paper. In what follows, we will denote by 
$\D{UC}(X)$, $\D{LSC}(X)$ and $\D{USC}(X)$ the space of uniformly continuous, lower semicontinuous and upper semicontinuous real functions on a metric space $X$, respectively.  

We are interested in a viscous Hamilton-Jacobi equation of the form
\begin{equation}\label{appB eq parabolic HJ}
\partial_{t }u=a(x) \partial^2_{xx} u +F(\partial_x u,x), \quad(t,x)\in (0,+\infty)\times\R,
\end{equation}
where  $F\in\CC(\R\times\R)$, and $a:\R\to [0,1]$ is a Lipschitz function satisfying the following assumption, for some constant $\hat\kappa_a > 0$: 
\begin{itemize}
\item[(A)]  $\sqrt{a}:\R\to [0,1]$\ is $\hat\kappa_a$--Lipschitz continuous.
\end{itemize}
We start with a comparison principle stated in a form which is the one we need in the paper.

\begin{prop}\label{appB prop comparison}
Suppose $a$ satisfies (A) and  $F\in\D{UC}\left(B_r\times\R\right)$ for every $r>0$. Let $v\in\D{USC}([0,T]\times\R)$ and $w\in\D{LSC}([0,T]\times\R)$ be, respectively, a sub and a supersolution of \eqref{appB eq parabolic HJ} in $(0,T)\times \R$ such that 
\begin{equation}\label{hyp 2}
\limsup_{|x|\to +\infty}\ \sup_{t\in [0,T]}\frac{v(t,x)-\theta x}{1+|x|}\leqslant 0 
\leqslant 
\liminf_{|x|\to +\infty}\ \sup_{t\in [0,T]}\frac{w(t,x)-\theta x}{1+|x|}
\end{equation}
for some $\theta\in\R$. Let us furthermore assume that either $\partial_x v$ or $\partial_x w$ belongs to $L^\infty\left((0,T)\times \R\right)$. Then, 
\[
v(t,x)-w(t,x)\leqslant \sup_{\R}\big(v(0,\,\cdot\,) - w(0,\,\cdot\,)\big)\quad\hbox{for every  $(t,x)\in (0,T)\times \R$.} 
\]
\end{prop}

\begin{proof}
The functions $\tilde v(t,x):=v(t,x)-\theta x$ and $\tilde w(t,x):=w(t,x)-\theta x$ are, respectively, a subsolution and a supersolution of \eqref{appB eq parabolic HJ} in $(0,T)\times \R$ with $F(\,\cdot\,,\theta +\,\cdot\,)$ in place of $F$. The assertion follows by applying \cite[Proposition 1.4]{D19} to $\tilde v$ and $\tilde w$.
\end{proof}

Let us now assume that $F$ belongs to the class $\Fam=\Fam(\alpha_0,\alpha_1,\gamma)$ of Hamiltonians satisfying the 
following set of assumptions, for some fixed constants $\alpha_0,\alpha_1>0$ and $\gamma>1$:\smallskip
\begin{itemize}
\item[(F1)] $\alpha_0|p|^\gamma-1/\alpha_0\leqslant F(p,x)\leqslant\alpha_1(|p|^\gamma+1)$\quad for all ${p,x\in\R}$;\medskip
\item[(F2)] $|F(p,x)-F(q,x)|\leqslant\alpha_1\left(|p|+|q|+1\right)^{\gamma-1}|p-q|$\quad for all $p,q,x\in\R$;\medskip
\item[(F3)] $|F(p,x)-F(p,y)|\leqslant\alpha_1\left(|p|^\gamma+1\right)|x-y|$\quad for all $p,x,y\in\R$. \medskip
\end{itemize}

The following holds. 

\begin{theorem}\label{appB teo well posed}
Suppose $a$ satisfies (A) and $F\in\Fam$. Then, for every $g\in\D{UC}(\R)$, there exists a unique function $u\in \D{UC}(\ccyl)$ that solves the equation \eqref{appB eq parabolic HJ}
subject to the initial condition $u(0,\,\cdot\,)=g$ on $\R$. If $g\in W^{2,\infty}(\R)$, then $u$ is Lipschitz continuous in $\ccyl$ and satisfies
\begin{equation*}
\|\partial_t u\|_{L^\infty(\ccyl)}\leqslant \kappa
\quad\text{and}\quad
\|\partial_x u\|_{L^\infty(\ccyl)}\leqslant \kappa
\end{equation*}
for some constant $\kappa$ that depends only on $\|g'\|_{L^\infty(\R)}$, $\|g''\|_{L^\infty(\R)}$, $\hat\kappa_a, \alpha_0,\alpha_1$ and $\gamma$. Furthermore, the dependence of $\kappa$ on $\|g'\|_{L^\infty(\R)}$ and $\|g''\|_{L^\infty(\R)}$ is continuous. 
\end{theorem}

\begin{proof}
A proof of this result when the initial datum $g$ is furthermore assumed to be bounded is given in \cite[Theorem 3.2]{D19}, see also \cite[Proposition 3.5]{AT}. 
This is enough, since we can always reduce to this case by possibly picking a function $\tilde g\in W^{3,\infty}(\R)\cap\CC^\infty(\R)$ such that $\|g-\tilde g\|_{L^\infty(\R)}<1$ (for instance, by mollification) and by considering equation \eqref{appB eq parabolic HJ} with
$\tilde F(p,x):=a(x)(\tilde g)'' +F(p + (\tilde g)',x)$ in place of $F$, and initial datum $g-\tilde g$. 
\end{proof}

For every fixed $\theta\in\R$, we will denote by $u_\theta$ the unique solution of \eqref{appB eq parabolic HJ} in $\D{UC}(\ccyl)$ satisfying  $u_\theta(0,x)=\theta x$ for all $x\in\R$. 
Theorem \ref{appB teo well posed} tells us that $u_\theta$ is Lipschitz in $\ccyl$ and that its Lipschitz constant depends continuously on $\theta$. In particular,  for every $r>0$, there exists a constant $\kappa_r>0$, depending only on 
$r, \hat\kappa_a, \alpha_0,\alpha_1$ and $\gamma$, such that 
\begin{equation}\label{appB eq Lipschitz estimate}
\|\partial_t u_\theta\|_{L^\infty(\ccyl)}\leqslant \kappa_r
\quad\text{and}\quad
\|\partial_x u_\theta\|_{L^\infty(\ccyl)}\leqslant \kappa_r
\quad 
\hbox{for every $\theta\in B_r$.} 
\end{equation}
Let us define 
\begin{eqnarray}\label{appB eq infsup}
	\HF^L(F) (\theta):=\liminf_{t\to +\infty}\ \frac{u_\theta(t,0)}{t}\quad\text{and}
	\quad
	\HF^U(F) (\theta):=\limsup_{t\to +\infty}\ \frac{u_\theta(t,0)}{t}.
\end{eqnarray}
By definition, we have \ $\HF^L(F) (\theta)\leqslant \HF^U(F) (\theta)$\ for all $\theta\in\R$. Furthermore, the following holds.

\begin{prop}\label{appB prop HF}
Suppose $a$ satisfies (A) and $F\in\Fam$. Then:
\begin{itemize}
\item[(i)] the functions  $\HF^L(F),\,\HF^U(F)$ satisfy (F1);\medskip
\item[(ii)] for every $r>0$, there exists a $\kappa_r>0$,\footnote{The constant $\kappa_r>0$ will be chosen according to \eqref{appB eq Lipschitz estimate}.} depending only on $r, \hat\kappa_a, \alpha_0,\alpha_1$ and $\gamma$, such that 
\[
|\HF^*(F)(\theta_1)-\HF^*(\theta_2)|
\leqslant
2\alpha_1\left(r+\kappa_r+1\right)^{\gamma-1}  |\theta_1-\theta_2|
\quad
\hbox{for all $\theta_1,\theta_2\in B_r$},
\]
where $\HF^*$ stands for either $\HF^L$ or $\HF^U$.
\end{itemize}
\end{prop}

\begin{proof} (i) Let us set 
\[
\alpha(h):=\alpha_0|h|^\gamma-\dfrac{1}{\alpha_0}
\quad\text{and}\quad
\beta(h):=\alpha_1\big( |h|^\gamma + 1\big)
\quad\hbox{for all $h\geqslant 0$.}
\]
Fix $\theta\in\R$. It is easily seen that the functions 
$v_\theta(t,x):=\theta x + t\alpha(|\theta|)$ and $w_\theta(t,x):=\theta x + t\beta(|\theta|)$ are, respectively, a (classical) sub and supersolution to \eqref{appB eq parabolic HJ}. By the comparison principle stated in Proposition \ref{appB prop comparison}, 
\[
v_\theta(t,x)\leqslant u_\theta (t,x) \leqslant w_\theta(t,x)
\quad
\hbox{for all $(t,x)\in\ccyl$.}
\]
The assertion follows by taking the liminf and limsup of the above inequalities, divided by $t$ and evaluated at $x=0$, as $t\to +\infty$.  

(ii) Let us fix $r>0$ and let $\kappa_r$ be a constant depending only on $r, \hat\kappa_a, \alpha_0,\alpha_1$ and $\gamma$, chosen according to \eqref{appB eq Lipschitz estimate}. Pick 
$\theta_1,\,\theta_2\in B_r$ and set 
\[
u_i(t,x):=u_{\theta_i}(t,x)-\theta_i x,
\quad
\hbox{$(t,x)\in\ccyl$, $i\in\{1,2\}$.}
\]
Each function $u_i$ is a Lipschitz solution to \eqref{appB eq parabolic HJ} with $F(\,\cdot\,,\theta_i+\,\cdot\,)$ in place of $F$ and with initial condition $u_i(0,\,\cdot\,) = 0$. Let us set 
$$M:=2\alpha_1\left(r+\kappa_r+1\right)^{\gamma-1}|\theta_1-\theta_2|.$$  
From the fact that $F$ satisfies (F2), we 
see that the following inequalities hold in the viscosity sense: 
\[
-M+F(x,\theta_2+\partial_x u_1)
\leqslant
F(x,\theta_1+\partial_x u_1)
\leqslant
M+F(x,\theta_2+\partial_x u_1)
\quad
\hbox{for all $(t,x)\in\cyl$.}
\]
We infer that the functions $v(t,x):=u_1(t,x)-Mt$ and $w(t,x):=u_1(t,x)+Mt$ are, respectively, a Lipschitz sub and supersolution to \eqref{appB eq parabolic HJ} with $F(\,\cdot\,,\theta_2+\,\cdot\,)$ in place of $F$. Furthermore, $v(0,\,\cdot\,) = w(0,\,\cdot\,) = 0$. By the comparison principle stated in Proposition \ref{appB prop comparison}, we infer that $v\leqslant u_2 \leqslant w$ in $\ccyl$, i.e.,
\[
u_{\theta_1} (t,x)-\theta_1 x-Mt 
\leqslant 
u_{\theta_2} (t,x)-\theta_2 x
\leqslant 
u_{\theta_1} (t,x)-\theta_1 x+Mt 
\quad
\hbox{for all $(t,x)\in\ccyl$.}
\]
The assertion follows by taking the liminf and limsup of the above inequalities, divided by $t$ and evaluated at $x=0$, as $t\to +\infty$.  
\end{proof}

According to \cite[Theorem 3.1]{DK17}, equation \eqref{appB eq parabolic HJ} homogenizes if and only if the functions 
$u^\eps_\theta(t,x):=\eps u_\theta(t/\eps,x/\eps)$ converge, locally uniformly in $\ccyl$, to a continuous function $\overline u_\theta(t,x)$. In this instance, we have $\HF^L(F)(\theta)=\HF^U(F)(\theta)$ and 
\[
\overline u_\theta(t,x)=\theta x+t\HF(F)(\theta),
\quad
\hbox{$(t,x)\in\ccyl$,}
\]
with $\HF(F)(\theta):=\HF^L(F)(\theta)=\HF^U(F)(\theta)$.

Here is a stability result for homogenization. 

\begin{theorem}\label{appB teo stability}
Suppose $a$ satisfies (A), the Hamiltonians $F$, $(F_n)_{n\in\N}$ belong to $\Fam$, and
\[
\lim_{n\to +\infty}\|F_n-F\|_{L^\infty(B_r\times\R)}=0
\quad
\hbox{for every $r>0$.}
\]
If equation \eqref{appB eq parabolic HJ} homogenizes for each $F_n$ with effective Hamiltonian $\HF(F_n)$, then it homogenizes for $F$, too, with effective Hamiltonian 
\[
\HF(F)(\theta)=\lim_{n\to +\infty} \HF(F_n)(\theta)
\quad
\hbox{for all $\theta\in\R$.}
\]
Furthermore, this convergence is locally uniform on $\R$. 
\end{theorem}

\begin{proof}
By hypothesis, we have $\HF(F_n):=\HF^L(F_n)=\HF^U(F_n)$ for each $n\in\N$. According to Proposition \ref{appB prop HF}, the effective Hamiltonians $\HF(F_n)$ make up a sequence of locally equi-bounded and locally equi-Lipschitz functions on $\R$, which is therefore pre-compact in $\CC(\R)$. Let us denote by $\overline F$ the limit of a convergent subsequence extracted from $\big(\HF(F_n)\big)_n$, which, for the sake of notational simplicity, we shall not relabel. We claim that, for every fixed $\theta\in\R$, 
\begin{equation}\label{appB eq main claim}
u^\eps_\theta(t,x):=\eps u_\theta(t/\eps,x/\eps) 
\underset{\eps\to 0^+}{\ucv}
\theta x +t\overline F(\theta)
\quad
\hbox{on $\ccyl$.}
\end{equation}
This will give us the desired result with $\HF(F)=\overline F$, yielding in particular that the whole sequence $\big(\HF(F_n)\big)_n$ converges to $\ol F$ in $\CC(\R)$. Let us prove the claim. Fix $\theta\in\R$ and, for each $n\in\N$, denote by $u^n_\theta$ the unique solution of 
\eqref{appB eq parabolic HJ} with $F_n$ in place of $F$ and satisfying the initial condition $u^n_\theta(0,x)=\theta x,\ x\in\R$. By \eqref{appB eq Lipschitz estimate} and our remark in the paragraph above it, there exists a common Lipschitz constant $\kappa>0$ for the functions $u^n_\theta$ on $\ccyl$. This implies that the functions 
\[
(u^n_\theta)^\eps(t,x):= \eps u^n_\theta(t/\eps,x/\eps),\quad (t,x)\in\ccyl,
\] 
are $\kappa$-Lipschitz in $\ccyl$, for every $n\in\N$ and $\eps>0$. Note that each $(u^n_\theta)^\eps$ is the unique solution in $\D{UC}(\ccyl)$ of
\begin{equation}\label{appB eq parabolic HJ n}
\partial_{t }u=\eps a\left(\frac x\eps\right) \partial^2_{xx} u +F_n\left(\partial_x u,\frac x\eps\right),\quad(t,x)\in(0,+\infty)\times\R,
\end{equation}
satisfying the initial condition $(u^n_\theta)^\eps(0,x)=\theta x,\ x\in\R$. Fix $r>\kappa$ and set 
$\delta_n:=\|F_n-F\|_{L^\infty(B_r\times\R)}$. The following inequalities hold in the viscosity sense:
\[
-\delta_n
\leqslant
F\left(\partial_x(u^n_\theta)^\eps, \frac x\eps\right)
-
F_n\left(\partial_x(u^n_\theta)^\eps, \frac x\eps\right)
\leqslant
\delta_n,
\quad
(t,x)\in\cyl.
\]
This readily implies that the functions $v(t,x):=(u^n_\theta)^\eps(t,x)-t\delta_n$ and $w(t,x):=(u^n_\theta)^\eps(t,x)+t\delta_n$ are, respectively, a sub and supersolution to \eqref{appB eq parabolic HJ n} with $F$ in place of $F_n$. Furthermore, they satisfy the initial condition $v(0,x)=w(0,x)=\theta x,\ x\in\R$. By the comparison principle, we infer that 
\[
(u^n_\theta)^\eps(t,x)-t\delta_n
\leqslant
u^\eps_\theta(t,x)
\leqslant
(u^n_\theta)^\eps(t,x)+t\delta_n
\quad
\hbox{for all $(t,x)\in\ccyl$.}
\]
By sending $\eps\to 0^+$, we get 
\[
\theta x +\left(\HF(F_n)(\theta)-\delta_n\right)t
\leqslant
\liminf_{\eps\to 0^+} u^\eps_\theta(t,x)
\leqslant
\limsup_{\eps\to 0^+} u^\eps_\theta(t,x)
\leqslant
\theta x +\left(\HF(F_n)(\theta)+\delta_n\right)t.
\]
Now we send $n\to +\infty$ to deduce that
\[
\liminf_{\eps\to 0^+} u^\eps_\theta(t,x)
=
\limsup_{\eps\to 0^+} u^\eps_\theta(t,x)
=
\theta x +t\overline F(\theta)
\quad
\hbox{for all $(t,x)\in\ccyl$}.
\]
This implies claim \eqref{appB eq main claim} since the functions $\{u^\eps_\theta:\eps>0\}$ are equi-Lipschitz on $\ccyl$. 
\end{proof}

\bibliographystyle{alpha}

\end{document}